\documentclass[a4paper]{amsart}
\usepackage{amsmath,amssymb,amsfonts}
\usepackage{exscale,url,doi,graphicx,xcolor,ulem}
\usepackage[linesnumbered,ruled,vlined,algosection]{algorithm2e}

\definecolor{color1}{rgb}{0,0,1}
\hypersetup{pdfborder={1 0 0}}

\def\diag{\ensuremath{\mathrm{diag}}}
\def\la{\ensuremath{\lambda}}
\def\RR{\ensuremath{\mathbb R}}
\def\CC{\ensuremath{\mathbb C}}
\def\spanl{\ensuremath{\mathrm{span}}}

\newtheorem{theorem}{Theorem}[section]
\newtheorem{lemma}{Lemma}[section]
\theoremstyle{definition}
\newtheorem{remark}{Remark}[section]
\newtheorem{definition}{Definition}[section]
\newtheorem{example}{Example}[section]
\numberwithin{equation}{section}

\SetKwInput{KwData}{Input}
\SetKwInput{KwResult}{Output}

\begin{document}

\title[Convergence analysis of BPSD-id]
      {Convergence analysis of a block preconditioned\\
       steepest descent eigensolver with implicit deflation}

\author[M. Zhou]{Ming Zhou}
\address{Universit\"at Rostock, Institut f\"ur Mathematik, 
        Ulmenstra{\ss}e 69, 18055 Rostock, Germany}
\email{ming.zhou@uni-rostock.de}

\author[Z. Bai]{Zhaojun Bai}
\address{Department of Computer Science and Department of Mathematics,
         University of California, Davis, CA 95616, USA}
\email{bai@cs.ucdavis.edu}

\author[Y. Cai]{Yunfeng Cai}
\address{Cognitive Computing Lab, Baidu Research,
         No. 10 Xibeiwang East Road, Beijing 100193, China}
\email{caiyunfeng@baidu.com}

\author[K. Neymeyr]{Klaus Neymeyr}
\address{Universit\"at Rostock, Institut f\"ur Mathematik, 
        Ulmenstra{\ss}e 69, 18055 Rostock, Germany}
\email{klaus.neymeyr@uni-rostock.de}


\keywords{
 Rayleigh quotient, gradient iterations, block eigensolvers.
 \hfill September 7, 2022}

\begin{abstract}
Gradient-type iterative methods for solving Hermitian eigenvalue problems
can be accelerated by using preconditioning and deflation techniques. 
A preconditioned steepest descent iteration with implicit deflation (PSD-id)
is one of such methods. 
The convergence behavior of the PSD-id is recently investigated 
based on the pioneering work of Samokish on the preconditioned 
steepest descent method (PSD).  The resulting non-asymptotic estimates indicate 
a superlinear convergence of the PSD-id under strong assumptions on the 
initial guess. The present paper utilizes an alternative convergence 
analysis of the PSD by Neymeyr under much weaker assumptions. 
We embed Neymeyr's approach into the analysis of the PSD-id using a 
restricted formulation of the PSD-id. More importantly, we extend the 
new convergence analysis of the PSD-id to a practically preferred 
block version of the PSD-id, or BPSD-id, and show the cluster robustness 
of the BPSD-id. Numerical examples are provided to validate the 
theoretical estimates. 
\end{abstract}

\maketitle

\section{Introduction} \label{sec:intro}

Determining the smallest eigenvalues and the associated eigenfunctions
of a self-adjoint elliptic partial differential operator is a common task 
in many application areas. The computational costs are usually high due to
large dimensions of discretized problems. By considering
an equivalent minimization problem for the Rayleigh quotient, one can develop 
suitable eigensolvers based on gradient iterations. The performance 
can be significantly improved by preconditioning techniques for 
better descent directions \cite{SAM1958},
and by a blockwise implementation 
for computing clustered eigenvalues \cite{BPK1996,OVT2006}. 
A further acceleration is enabled by modifying the preconditioners 
with certain shifts after deflation.
These methodological improvements are particularly meaningful
for ill-conditioned eigenvalue problems
arising from applications, such as the partition-of-unity finite element method 
for solving the Kohn-Sham equation in electronic structure calculations,
where the target eigenvalues are not well separated from the rest of 
the spectrum; see \cite{CBPS2018} and references therein.

In this paper, we consider the generalized Hermitian definite matrix eigenvalue problem
\begin{equation}\label{evp}
Hu=\lambda Su
\end{equation}
where $H,S\in\CC^{n \times n}$ are Hermitian, and $S$ is positive definite.
The eigenvalues of $(H,S)$ are all real. We denote them by
$\la_1,\ldots,\la_n$ and arrange them as $\la_1\le\cdots\le\la_n$.
We aim at computing
the first $m$ eigenvalues together with the associated eigenvectors.
Typically, $m \ll n$. 

When $m=1$, the task can be restated as minimizing
the Rayleigh quotient
\begin{equation}\label{rq}
\rho: \CC^n{\setminus}\{0\}\to\RR,\quad \rho(z)=\frac{z^*Hz}{z^*Sz}
\end{equation}
where the superscript ${}^*$ denotes conjugate transpose.
The preconditioned steepest descent iteration (PSD) 
\begin{equation}\label{psd}
z^{(\ell+1)}=z^{(\ell)}-\omega^{(\ell)}Kr^{(\ell)}
\end{equation}
is applicable to the minimization problem of \eqref{rq}.
Therein $K\in\CC^{n \times n}$ is a Hermitian preconditioner,
and the residual $r^{(\ell)}=Hz^{(\ell)}-\rho(z^{(\ell)})Sz^{(\ell)}$
is collinear with the gradient vector $\nabla\rho(z^{(\ell)})$
of the Rayleigh quotient $\rho(z^{(\ell)})$.
The descent direction $-Kr^{(\ell)}$ generalizes $-\nabla\rho(z^{(\ell)})$,
and is expected to enable an acceleration.
An optimal stepsize $\omega^{(\ell)}$ can be determined
by the Rayleigh-Ritz procedure applied on the subspace
$\spanl\{z^{(\ell)},Kr^{(\ell)}\}$. 

The choice of $K$ significantly affects
the convergence rate of the PSD. 
Setting $K$ equal to the $n{\times}n$ identity matrix $I$
leads to a slow and mesh-depedent convergence 
for discretized operator eigenvalue
problems \cite{KNN2003e}. A shift-and-invert preconditioner
$K=(H-\sigma S)^{-1}$ with $\sigma<\la_1$ is much more efficient.
According to the analyses of gradient eigensolvers 
in \cite{KNS1991,NOZ2011}, one can derive the single-step convergence factor
\[
 \kappa=\left(\frac{\la_1-\sigma}{\la_2-\sigma}\right)
 \left(\frac{\lambda_n-\lambda_2}{\lambda_n-\lambda_1}\right)
\]
for approximate eigenvectors (with respect to tangent values of error angles)
and additionally $\kappa^2/(2-\kappa)^2$ for approximate eigenvalues in the final phase
(with respect to relative positions between $\la_1$ and $\la_2$).
These convergence factors are bounded away from $1$ when $\la_1\ll\la_2$ and
can be improved by empirically increasing $\sigma$.

In practice, the shift-and-invert preconditioner $K$ is 
often implemented approximately,
e.g., by iteratively solving linear systems of the form 
\begin{equation}\label{linsys}
 \widetilde{H}p=r
 \quad\mbox{with}\quad
 \widetilde{H}=H-\sigma S.
\end{equation}
Thus it is more important to analyze 
a preconditioner $K\approx\widetilde{H}^{-1}$. 
The analysis from \cite{NEY2012}
begins with a convergence argument for \eqref{linsys}
with respect to the error propagation matrix $I-K\widetilde{H}$.
By using a quality parameter $\varepsilon$ defined in the condition
\begin{equation}\label{cond}
 \|I-K\widetilde{H}\|_{\widetilde{H}}\le\varepsilon<1,
\end{equation}
the eigenvalue convergence factor $\kappa^2/(2-\kappa)^2$
for the special case $K=\widetilde{H}^{-1}$ (where $\varepsilon=0$)
is generalized as
\begin{equation}\label{cvf}
 \big(\kappa+\varepsilon(2-\kappa)\big)^2/\big((2-\kappa)+\varepsilon\kappa\big)^2.
\end{equation}
As we can see, the convergence rate of the PSD can be deteriorated by
$\varepsilon\approx 1$ or $\kappa \approx 1$.
One can reduce $\varepsilon$ by using a proper
linear system solver. However, reducing $\kappa$ is difficult
if $\la_1$ and $\la_2$ are tightly clustered.
In this case, the PSD needs to be modified, e.g., by a blockwise implementation
or certain subspace extensions. The resulting iterations also provide
approximations of further eigenvalues.

For $m>1$, computing the $m$ smallest eigenvalues of
the matrix pair $(H,S)$ corresponds to the trace minimization of $Z^*HZ$
among all $S$-orthonormal matrices $Z\in\CC^{n \times m}$ \cite{STS1990}.
A typical approach is the block version of the PSD, BPSD in short.
Therein each iterate $Z^{(\ell)}$ is an $S$-orthonormal Ritz basis matrix, 
and $\Theta^{(\ell)}=(Z^{(\ell)})^*HZ^{(\ell)}$ is a diagonal matrix
whose diagonal consists of the Ritz values in $\spanl\{Z^{(\ell)}\}$.
The next iterate $Z^{(\ell+1)}$ is determined by
applying the Rayleigh-Ritz procedure to the trial subspace
\[
 \spanl\{Z^{(\ell)},KR^{(\ell)}\}
 \quad\mbox{with}\quad
 R^{(\ell)}=HZ^{(\ell)}-SZ^{(\ell)}\Theta^{(\ell)}.
\]
The columns of $Z^{(\ell+1)}$ are
$S$-orthonormal Ritz vectors associated with the $m$ smallest Ritz values
in $\spanl\{Z^{(\ell)},KR^{(\ell)}\}$. The convergence behavior of the BPSD
with respect to individual Ritz values
has been analyzed in \cite{NEZ2014}.
The convergence factor for the $i$th Ritz value
in the final phase is given by \eqref{cvf} with the generalized parameter
\[
 \kappa=\left(\frac{\la_i-\sigma}{\la_{i+1}-\sigma}\right)
 \left(\frac{\lambda_n-\lambda_{i+1}}{\lambda_n-\lambda_i}\right).
\]

An alternative convergence factor presented in \cite{NEZ2019}
concerns the non-optimized version $Z^{(\ell+1)}$ $=Z^{(\ell)}-KR^{(\ell)}$
of the BPSD, and depends on the ratio $(\la_i-\sigma)/(\la_{m+1}-\sigma)$
instead of $(\la_i-\sigma)/(\la_{i+1}-\sigma)$.
This reflects the cluster robustness of the BPSD
for $\la_i \ll \la_{m+1}$ which can be ensured by enlarging the block size $m$.

Deflation is required for block iterative eigensolvers
if the number of target eigenvalues exceeds the block size. In the (B)PSD,
once the residual of an approximate eigenvector is sufficiently small,
one can store it in a basis matrix $U$ consisting of accepted
$S$-orthonormal approximate eigenvectors.
The further iterates are $S$-orthogonalized against $U$
so that they converge toward eigenvectors associated with
the next eigenvalues. Such an orthogonalization is called
implicit deflation \cite[Section 6.2.3]{SAA1992}.

Recently, a combination of the PSD with implicit deflation (PSD-id)
has been investigated in \cite{CBPS2018}. A remarkable feature
of the PSD-id is that the preconditioner is variable depending on the current
approximate eigenvalue, somewhat similarly to the Jacobi-Davidson method
or an inexact Rayleigh quotient iteration \cite{NTY2003}. This accelerates
the convergence in the final phase in comparison to a fixed preconditioner.

The convergence behavior of the PSD-id is analyzed in \cite{CBPS2018}
based on the pioneering work of Samokish on the PSD \cite{SAM1958}.
In particular, \cite[Theorem 3.2]{CBPS2018} extends
a non-asymptotic reformulation of the classical estimate \cite[(10)]{SAM1958}
presented by Ovtchinnikov in \cite[Theorem 2.1]{OVT2006s}.
Precisely, a convergence rate bound for approximating
an interior eigenvalue $\la_i$ $(>\la_{i-1})$ with the PSD-id
and a variable preconditioner $K$ is derived by using the matrix
$M=Q^*(H - \la_i S)Q$,
where $Q$ is an $S$-orthogonal projector
onto the invariant subspace associated with the eigenvalues
$\la_i,\ldots,\la_n$. The resulting estimate contains
an essential parameter
\[
 q=\frac{\beta-\alpha}{\beta+\alpha}, 
\]
where $\alpha$ and $\beta$ are the smallest positive and largest eigenvalues
of the matrix product $KM$. A corresponding asymptotic estimate
then uses $q^2$ as the convergence factor. It is possible to reduce
$q^2$ by modifying $K$ with proper shifts for a superlinear convergence.

A drawback of \cite[Theorem 3.2]{CBPS2018} is its technical assumption
on parameters associated with $K$ and $M$.
The assumption is rather restricted as
the current approximate eigenvalue has to be sufficiently close to $\la_i$.
More significantly, it would be extremely difficult, if possible, to generalize
the analysis to a block version of the PSD-id, or BPSD-id,
which is the algorithm used in practice.

In this paper, we first recast the PSD-id as the PSD applied to 
a restricted eigenvalue problem, and then utilize the approaches 
proposed in \cite{NEY2012} to provide a convergence analysis of the PSD-id.
More importantly, we are able to present a convergence analysis of the BPSD-id
by extending the convergence analysis of the BPSD in \cite{NEZ2014}. 
The resulting sharp single-step estimates of the PSD-id and the BPSD-id
are presented in Theorems~\ref{psdidth} and~\ref{bpsdidth}.
Further estimates in Sections~\ref{sec:sigma} and~\ref{sec:multi} 
deal with larger shifts and clustered eigenvalues.
Numerical examples are provided to verify the sharpness of 
convergence estimates. These theoretical results and numerical experiments
advance our in-depth understanding of the convergence behaviors of 
the PSD-id and the BPSD-id.

For ease of references, we restate the estimates 
of the PSD and the BPSD in \cite{NEY2012,NEZ2014} for 
a standard Hermitian eigenvalue problem with notation 
that are compatible with a restricted formulation
of the PSD-id and the BPSD-id in Section~\ref{sec:prelim}.

\begin{theorem}[\cite{NEY2012,NEZ2014}]\label{thmbpsd}
Consider a Hermitian positive definite matrix $A\in\CC^{d \times d}$
together with its eigenvalues $\zeta_1\le\cdots\le\zeta_d$ and 
a Hermitian positive definite preconditioner
$T\in\CC^{d \times d}$ satisfying
\begin{equation}\label{pcq}
 \|I-TA\|_{A}\le\varepsilon<1.
\end{equation}
If $\eta = (x^*Ax)/(x^*x)$ for a nonzero vector $x\in\CC^d$
is located in an eigenvalue interval
$(\zeta_j,\zeta_{j+1})$, then it holds for
the smallest Ritz value $\eta'$ of $A$
in $\spanl\big\{x,\,T(Ax-\eta x)\big\}$ that
\begin{equation}\label{estpsd}
 \frac{\eta'-\zeta_j}{\,\zeta_{j+1}-\eta'\,}
 \le\left(\frac{\kappa+\varepsilon(2-\kappa)}{(2-\kappa)+\varepsilon\kappa}\right)^2
 \frac{\eta-\zeta_j}{\,\zeta_{j+1}-\eta\,}
 \quad\mbox{with}\quad
 \kappa=\frac{\zeta_j(\zeta_d-\zeta_{j+1})}{\,\zeta_{j+1}(\zeta_d-\zeta_j)\,}.
\end{equation}
The equality in \eqref{estpsd} is attainable
in the limit case \,$\eta\to\zeta_j$\,
in an invariant subspace associated with
the eigenvalues $\zeta_j$, $\zeta_{j+1}$ and $\zeta_d$.

Consider further a $c$-dimensional subspace $\mathcal{X}\subset\CC^d$
and orthonormal Ritz vectors $x_1,\ldots,x_c$ associated with
the Ritz values $\eta_1\le\cdots\le\eta_c$ of $A$ induced by $\mathcal{X}$.
Denote by $\eta'_1\le\cdots\le\eta'_c$ the $c$ smallest Ritz values of $A$
induced by $\spanl\big\{X,\,T(AX-XD)\big\}$
with $X=[x_1,\ldots,x_c]$ and $D=\diag(\eta_1,\ldots,\eta_c)$.
For each $t\in\{1,\ldots,c\}$, if $\eta_t$ is located in an eigenvalue interval
$(\zeta_j,\zeta_{j+1})$, then it holds that
\begin{equation}\label{estbpsd}
 \frac{\eta'_t-\zeta_j}{\,\zeta_{j+1}-\eta'_t\,}
 \le\left(\frac{\kappa+\varepsilon(2-\kappa)}{(2-\kappa)+\varepsilon\kappa}\right)^2
 \frac{\eta_t-\zeta_j}{\,\zeta_{j+1}-\eta_t\,}
 \quad\mbox{with}\quad
 \kappa=\frac{\zeta_j(\zeta_d-\zeta_{j+1})}{\,\zeta_{j+1}(\zeta_d-\zeta_j)\,}.
\end{equation}
The equality in \eqref{estbpsd} is attainable
in the limit case \,$\eta_t\to\zeta_j$\,
in an invariant subspace associated with
the eigenvalues $\zeta_j$, $\zeta_{j+1}$ and $\zeta_d$.
\end{theorem}

The remaining part of this paper is organized as follows. 
In Section~\ref{sec:prelim}, we introduce
the algorithmic structure of the (B)PSD-id, and present 
a restricted formulation
as the starting point of our analysis. 
In Section~\ref{sec:single}, the estimates of 
the (B)PSD in Theorem \ref{thmbpsd} are 
applied to certain representations
of the (B)PSD-id in an invariant subspace 
in order to derive sharp single-step estimates.
Multi-step estimates are presented
in Section~\ref{sec:multi} for the cluster robustness of the BPSD-id.
The theoretical convergence estimates of the BPSD-id 
are demonstrated by numerical experiments in Section~\ref{sec:numerics}.
 
\section{PSD-id and BPSD-id algorithms and restricted formulations} 
\label{sec:prelim}

\subsection{PSD-id and BPSD-id algorithms} \label{sec:algs}

In this section, we reformulate the original algorithm of 
the PSD-id from \cite{CBPS2018} in Algorithm \ref{psdid}.
Therein we drop the step indices of the iterates since they are not
needed in the derivation of our new estimates.
A few remarks of Algorithm \ref{psdid} are in order.

\begin{itemize} 
\item Line 1:
The $S$-orthogonalization against the $i-1$ accepted approximate eigenvectors
is made once at initialization and then automatically within the
computation of Ritz pairs.

\item Line 2:
The stopping criterion uses the norm $\|\cdot\|_{S^{-1}}$. It 
implies that there exists an eigenvalue $\lambda$ of $(H,S)$ which 
fulfills
\[
 |\lambda-\rho(z)|\le\|r\|_{S^{-1}}/\|z\|_S,
\]
cf.~\cite[Theorem 15.9.1]{PAR1980}. 

\item Line 3: The preconditioner $K$ is constructed to be
effectively positive definite for ensuring the convergence; see Definition \ref{epd}.
In practice, $p=Kr$ can be implemented by approximate solution
of the linear system $(H-\sigma S)p=r$.
The shift $\sigma$ will be discussed in the convergence analysis
and the numerical experiments.
In principle we set $\sigma$ slightly smaller than
the $i$th smallest eigenvalue $\la_i$,
or equal to the current approximate eigenvalue $\rho(z)$
if it is close to $\la_i$.
\end{itemize} 

\begin{algorithm}
\caption{$z=\mbox{PSD-id}\,(U,z)$}\label{psdid}
\KwData{$S$-orthonormal matrix $U\in\CC^{n\times(i-1)}$ whose columns are
 accepted approximate eigenvectors associated 
 with the $i-1$ smallest eigenvalues;
 initial guess $z\in\CC^n$.}
\KwResult{approximate eigenvector $z$
associated with the $i$th smallest eigenvalue.}
$S$-orthogonalize $z$ against $U$: $z = z - U U^* S z$;
\quad $z=z/\|z\|_S$; 
\quad $r=Hz-\rho(z)Sz$\;
\While{$\|r\|_{S^{-1}}$ not sufficiently small}{
compute a preconditioned residual $p=Kr$\;
update $z$ by an $S$-normalized Ritz vector
associated with the $i$th smallest Ritz value
in the subspace $\spanl\{U,z,p\}$\;
$r=Hz-\rho(z)Sz$\;}
\end{algorithm}

Algorithm \ref{bpsdid} describes the BPSD-id. Therein
the usage of a suitable block size $\widetilde{k}$ can overcome the possible
convergence stagnation of the PSD-id for clustered eigenvalues;
see Section~\ref{sec:multi}.
A few remarks of Algorithm \ref{bpsdid} are in order.

\begin{itemize} 
\item Line 4: The block residual $R$
actually consists of the residuals of individual Ritz vectors.
The first $k$ columns of $R$ are considered in the stopping criterion.

\item There are two implementations with different choices of the block size
for computing the smallest eigenvalues of $(H,S)$
and the associated eigenvectors, beginning with
a random initial guess $W\in\CC^{n \times m}$
where $m$ is larger than the number of target eigenvalues.
A straightforward implementation
with a fixed block size $\widetilde{k}<m$ has the form 
\[
 \mbox{BPSD-id}\,
 \big(W(:,\,1\,{:}\,i{-}1),W(:,\,i\,{:}\,i{-}1{+}\widetilde{k}),k\big),
\]
i.e., each outer step only treats $\widetilde{k}$ columns of $W$.
The leading index $i$ of these columns is initially $1$
and will be updated together with $W(:,\,1\,{:}\,i{-}1)$
by using already converged columns.
An alternative implementation uses the block size
$\widetilde{k}=m{-}i{+}1$ depending on the current index $i$.
Thus $W$ is entirely modified as early as in the first outer step.
The latter columns of $W$ can provide more accurate initial data
for the inner loop in the next outer step.
\end{itemize} 

\begin{algorithm}
\caption{$Z=\mbox{BPSD-id}\,(U,Z,k)$}\label{bpsdid}
\KwData{$S$-orthonormal matrix $U\in\CC^{n\times(i-1)}$ whose columns are
 accepted approximate eigenvectors associated with the $i-1$ 
 smallest eigenvalues;
 initial guess $Z\in\CC^{n \times \widetilde{k}}$,
 $\widetilde{k}$ is the block size.}
\KwResult{matrix $Z\in\CC^{n \times k}$
($k\le\widetilde{k}$) consisting of approximate eigenvectors
associated with the $i$th to the $(i{-}1{+}k)$th smallest eigenvalues.}
$S$-orthogonalize $Z$ against $U$\;
update $Z$ by $S$-orthonormal Ritz vectors in $\spanl\{Z\}$\;
$R=HZ-SZ(Z^*HZ)$\;
\While{$\|R(:,\,1\,{:}\,k)\|_{S^{-1}}$ not sufficiently small}{
compute a block preconditioned residual $P=KR$\;
update $Z$ by $S$-orthonormal Ritz vectors associated with
the $i$th to the $(i{-}1{+}\widetilde{k})$th smallest Ritz values
in the subspace $\spanl\{U,Z,P\}$\;
$R=HZ-SZ(Z^*HZ)$\;}
$Z=Z(:,\,1\,{:}\,k)$\;
\end{algorithm}

\subsection{Restricted formulations of the PSD-id and the BPSD-id} 
\label{sec:reform}

As in \cite{CBPS2018},
we assume by ignoring sufficiently small numerical errors
that the columns $u_1$, $\ldots$ , $u_{i-1}$ of the matrix $U$
in the PSD-id (Algorithm \ref{psdid})
and the BPSD-id (Algorithm \ref{bpsdid}) are
exact $S$-orthonormal eigenvectors associated with
the $i{\,-\,}1$ smallest eigenvalues $\lambda_1\le\cdots\le\lambda_{i-1}$
of $(H,S)$,  and the iterate $z$ or $Z$ after the first $S$-orthogonalization
against $U$ has full rank.

As the starting point of our new convergence analysis,
we represent the PSD-id by the PSD applied to 
a restricted eigenvalue problem.
Let us first extend $U$ by $V=[u_i,\ldots,u_n]$
as an $S$-orthonormal basis of $\CC^n$
where $u_i,\ldots,u_n$ are eigenvectors associated with
the remaining eigenvalues $\lambda_i\le\cdots\le\lambda_n$, i.e.,
\begin{equation}\label{bevp}
HV=SV\Lambda 
\quad\mbox{with}\quad
\Lambda=\mathrm{diag}(\lambda_i,\ldots,\lambda_n)
\quad\mbox{and}\quad \spanl\{V\}=\spanl\{U\}^{\perp_S}.
\end{equation}

For the PSD-id, 
an arbitrary $\widetilde{z}\in\spanl\{V\}{\setminus}\{0\}$ can be represented by
\begin{equation}\label{vecr}
 \widetilde{z}=V\widetilde{c}
 \quad\mbox{with}\quad
 \widetilde{c}=V^*S\widetilde{z}.
\end{equation}
By using the identity \eqref{bevp}, 
a relation between the Rayleigh quotient $\rho(\cdot)$ 
of $(H, S)$ defined in \eqref{rq} and 
the restricted Rayleigh quotient of $\Lambda$: 
\begin{equation}\label{lrq}
\widetilde{\rho}: \CC^{n-i+1}{\setminus}\{0\}\to\RR,\quad
\widetilde{\rho}(w)=\frac{w^*\Lambda w}{w^*w}
\end{equation}
is given by
\begin{equation}\label{rqe}
\rho(\widetilde{z})=\rho(V\widetilde{c})
=\frac{\widetilde{c}^*V^*HV\widetilde{c}}{\widetilde{c}^*V^*SV\widetilde{c}}
\stackrel{\eqref{bevp}}{=}\frac{\widetilde{c}^*V^*SV\Lambda\widetilde{c}}{\widetilde{c}^*\widetilde{c}}
=\frac{\widetilde{c}^*\Lambda\widetilde{c}}{\widetilde{c}^*\widetilde{c}}=\widetilde{\rho}(\widetilde{c}).
\end{equation}
Consequently, the target eigenvalue of the PSD-id (Algorithm \ref{psdid})
can be interpreted by
\[
 \min_{\widetilde{z}\in\CC^n{\setminus}\{0\},\ 
  U^* S \widetilde{z} = 0} \rho(\widetilde{z})
 = \min_{\widetilde{z} \in \spanl\{V\}{\setminus}\{0\}} \rho(\widetilde{z}) 
 = \min_{\widetilde{c} \in \mathbb{C}^{n-i+1}{\setminus}\{0\}}\widetilde{\rho}(\widetilde{c}).
\] 
It implies that 
the PSD-id for computing the $i$th smallest eigenvalue of $(H,S)$
is equivalent to the PSD for computing the smallest eigenvalue of $\Lambda$.
The following Lemma presents such relationship in detail.

\begin{lemma}\label{lmpsdid}
Denote by $z$ and $z'$ two successive iterates of the PSD-id.
\begin{itemize}
\setlength\itemsep{1ex}
\item[(i)] If $z$ belongs to $\spanl\{V\}{\setminus}\{0\}$, then also $z'$.
\item[(ii)] Let $c=V^*Sz$ and $c'=V^*Sz'$
be the coefficient vectors of $z$ and $z'$
with respect to the representation \eqref{vecr}.
Then $c'$ is a minimizer of \,$\widetilde{\rho}(\cdot)$ in
$\spanl\{c,\widetilde{K}\widetilde{r}\}$ 
with $\widetilde{K}=V^*SKSV$ and 
$\widetilde{r}=\Lambda c-\widetilde{\rho}(c)c$.
\end{itemize} 
\end{lemma}
\begin{proof}
(i) 
If $z$ belongs to $\spanl\{V\}{\setminus}\{0\}$,
then the dimension of the trial subspace $\spanl\{U,z,p\}$ is at least $i$.
This verifies the existence of the $i$th smallest Ritz value
in $\spanl\{U,z,p\}$ and the existence of
the next iterate $z'$ which is an associated Ritz vector.
Moreover, the columns of $U$ are eigenvectors
associated with the $i-1$ smallest eigenvalues
and automatically Ritz vectors associated with 
the $i-1$ smallest Ritz values in $\spanl\{U,z,p\}$.
Thus $z'$ is $S$-orthogonal to $\spanl\{U\}$
and belongs to $\spanl\{V\}{\setminus}\{0\}$.

(ii) 
By using the representation $z=Vc$
and the projector $Q=VV^*S$ onto $\spanl\{V\}$, we get
\[
\spanl\{U,z,p\}=\spanl\{U\}\oplus\spanl\{Vc,Qp\}
\]
and
\[
\begin{split}
Qp \,\ = \,\ & (VV^*S)Kr = VV^*SK\big(Hz-\rho(z)Sz\big)
\stackrel{\eqref{rqe}}{=}VV^*SK\big(HVc-\widetilde{\rho}(c)SVc\big)\\[1ex]
\stackrel{\eqref{bevp}}{=}\,& VV^*SK\big(SV\Lambda c-\widetilde{\rho}(c)SVc\big)
=VV^*SKSV\big(\Lambda c-\widetilde{\rho}(c)c\big)
=V\widetilde{K}\widetilde{r}. 
\end{split}
\]
Therefore
\[
 \spanl\{U,z,p\}
 =\spanl\{U\}\oplus\spanl\{Vc,V\widetilde{K}\widetilde{r}\}
 =\spanl\{U\}\oplus V \cdot \spanl\{c,\widetilde{K}\widetilde{r}\}.
\]
Recall that the columns of $U$ are Ritz vectors associated with 
the $i-1$ smallest Ritz values in $\spanl\{U,z,p\}$,
the $i$th smallest Ritz value is just the smallest Ritz value
in $V \cdot \spanl\{c,\widetilde{K}\widetilde{r}\}$.
The associated Ritz vector $z'$ is thus
a minimizer of $\rho(\cdot)$ therein.
The relation \eqref{rqe} ensures that minimizing $\rho(\cdot)$
in the subspace $V \cdot \spanl\{c,\widetilde{K}\widetilde{r}\}$
is equivalent to minimizing $\widetilde{\rho}(\cdot)$
in the ``coefficient subspace'' $\spanl\{c,\widetilde{K}\widetilde{r}\}$.
Consequently,  the coefficient vector $c'$ of $z'$ is a minimizer
of $\widetilde{\rho}(\cdot)$ in $\spanl\{c,\widetilde{K}\widetilde{r}\}$.
\end{proof}

Lemma \ref{lmpsdid} indicates that
$c'$ is a Ritz vector associated with the smallest Ritz value
of $\Lambda$ in $\spanl\{c,\widetilde{K}\widetilde{r}\}$.
Thus $c$ and $c'$ can be regarded as two successive iterates
of a PSD iteration for minimizing $\widetilde{\rho}(\cdot)$
with $\widetilde{K}$ as preconditioner.
Consequently, we can analyze the PSD-id in terms of the PSD iteration 
with successive iterates $c$ and $c'$. 

Now let us represent the BPSD-id by the BPSD applied to 
a restricted eigenvalue problem.
For an arbitrary matrix $\widetilde{Z}\in\CC^{n \times l}$ having full rank and
satisfying $\spanl\{\widetilde{Z}\}\subseteq\spanl\{V\}$, one can define
$\widetilde{C}=V^*S\widetilde{Z}$ so that $\widetilde{Z}=V\widetilde{C}$.
Based on the relations
\begin{equation}\label{rve}
 \widetilde{Z}^*H\widetilde{Z}=\widetilde{C}^*V^*HV\widetilde{C}
 \stackrel{\eqref{bevp}}{=}\widetilde{C}^*V^*SV\Lambda\widetilde{C}
 =\widetilde{C}^*\Lambda\widetilde{C}
 \quad\mbox{and}\quad
 \widetilde{Z}^*S\widetilde{Z}
 =\widetilde{C}^*V^*SV\widetilde{C}=\widetilde{C}^*\widetilde{C},
\end{equation}
the Ritz values of $(H,S)$ in $\spanl\{\widetilde{Z}\}$
coincide with those of $\Lambda$ in $\spanl\{\widetilde{C}\}$.
The respective Ritz vectors
can be converted by multiplications with $V$ or $V^*S$
analogously to \eqref{vecr}.

\begin{lemma}\label{lmbpsdid}
Denote by $Z$ and $Z'$ two successive iterates of the BPSD-id.
\begin{itemize}
\setlength\itemsep{1ex}
\item[(i)] If $Z$ has full rank and all its $\widetilde{k}$ columns
belong to $\spanl\{V\}$, then also $Z'$.

\item[(ii)] Define the coefficient matrices
$C=V^*SZ$ and $C'=V^*SZ'$.
Then the columns of \,$C$ are orthonormal Ritz vectors
of $\Lambda$ in $\spanl\{C\}$, i.e.,
$\Theta=C^*\Lambda C$ is a diagonal matrix
whose diagonal entries are the corresponding Ritz values.
Furthermore, the columns of \,$C'$
are orthonormal Ritz vectors associated with
the $\widetilde{k}$ smallest Ritz values of $\Lambda$ in
$\spanl\{C,\widetilde{K}\widetilde{R}\}$ for
$\widetilde{K}=V^*SKSV$ and
$\widetilde{R}=\Lambda C - C\Theta$.
\end{itemize}
\end{lemma}
\begin{proof}
(i)
The dimension of $\spanl\{U,Z,P\}$ is at least
$i{\,-\,}1{\,+\,}\widetilde{k}$ so that $Z'$ is constructed
by $S$-orthonormal Ritz vectors associated with
the $i$th to the $(i{-}1{+}\widetilde{k})$th smallest Ritz values
of $(H,S)$ in $\spanl\{U,Z,P\}$, and has full rank.
These Ritz vectors (columns of $Z'$) are $S$-orthogonal
to those associated with the $i-1$ smallest eigenvalues
and thus belong to the $S$-orthogonal complement of $\spanl\{U\}$,
i.e., $\spanl\{V\}$.

(ii)
The statement for $C$ is simply based on \eqref{rve}
applied to $Z$ and $C$ together with the fact that the columns of $Z$
are $S$-orthonormal Ritz vectors of $(H,S)$ in $\spanl\{Z\}$.
In order to verify the statement for $C'$, the relation
\[
 \spanl\{U,Z,P\}=\spanl\{U\}\oplus\spanl\{Z,QP\}
\]
with the projector $Q=VV^*S$ onto $\spanl\{V\}$
indicates that the columns of $Z'$ are
$S$-orthonormal Ritz vectors associated with
the $\widetilde{k}$ smallest Ritz values of $(H,S)$
in $\spanl\{Z,QP\}$. Moreover, by using
$\Theta=C^*\Lambda C=Z^*HZ$, it holds that
\[
\begin{split}
QP & = (VV^*S)KR=VV^*SK\big(HZ-SZ\Theta\big)
=VV^*SK\big(HVC-SVC\Theta\big)\\
 & = VV^*SK\big(SV\Lambda C-SVC\Theta\big)
=VV^*SKSV\big(\Lambda C-C\Theta\big)
=V\widetilde{K}\widetilde{R}.
\end{split}
\]
so that $\spanl\{Z,QP\}=V \cdot \spanl\{C,\widetilde{K}\widetilde{R}\}$.
Applying \eqref{rve} to $\spanl\{Z,QP\}$
and $\spanl\{C,\widetilde{K}\widetilde{R}\}$
completes the verification.
\end{proof}

By Lemma~\ref{lmbpsdid}, 
we can analyze the convergence behavior of the BPSD-id
in terms of the BPSD iteration with successive iterates $C$ and $C'$.

\section{Sharp single-step estimates} \label{sec:single} 

In this section, 
we first analyze the convergence behavior
of the PSD-id and the BPSD-id.
Section~\ref{sec:psdid} presents an alternative convergence analysis
of the PSD-id (Algorithm \ref{psdid}) in comparison to \cite{CBPS2018}.
The estimate \eqref{estpsd} of the PSD with weaker assumptions is applied
to the restricted formulation of the PSD-id introduced in Lemma \ref{lmpsdid}.
This results in a sharp single-step estimate of the PSD-id. In Section \ref{sec:bpsdid},
the convergence of the BPSD-id (Algorithm \ref{bpsdid}) is investigated
by using the estimate \eqref{estbpsd}
together with Lemma \ref{lmbpsdid}. 
In Section \ref{sec:sigma}, we discuss an extension of
the main results under the notion of so-called larger shifts.

In preparation for the main analysis in this section, 
we characterize the preconditioner $K$
for the PSD-id and the BPSD-id with respect to its restricted form
$\widetilde{K}$ arising from Lemmas \ref{lmpsdid} and \ref{lmbpsdid}.

\begin{definition}[\mbox{\cite[Defintion 2.1]{CBPS2018}}]\label{epd}
A preconditioner $K$ is called {\em effectively positive definite},
if $\widetilde{K}=V^*SKSV\in\CC^{(n-i+1)\times(n-i+1)}$ 
is positive definite, where $V$ is defined in \eqref{bevp}.  
$\widetilde{K}$ is called an {\em effective form} of $K$. 
\end{definition}

A typical example of effectively positive definite preconditioners
is the shift-invert preconditioner $K=(H-\sigma S)^{-1}$
where the shift $\sigma$ is smaller than $\la_i$,
and not an eigenvalue of $(H,S)$. In this case, the corresponding effective form
$\widetilde{K}$ is actually a diagonal matrix since
by~\eqref{bevp}, 
\begin{align} \label{irf}  
\widetilde{K}  = V^*SKSV=V^*S(H-\sigma S)^{-1}SV
 =V^*SV(\Lambda - \sigma \widetilde{I})^{-1}
=(\Lambda - \sigma \widetilde{I})^{-1}, 
\end{align} 
where $\widetilde{I}=I_{n-i+1}$. $\widetilde{K}$ is positive definite
since $(\la_j-\sigma)^{-1}>0$ for each $j \ge i$ due to $\sigma<\la_i$.

\subsection{Sharp single-step estimate of the PSD-id} \label{sec:psdid}

We first provide a simple proof on the monotonicity of
the approximate eigenvalues, which has been proven
in a cumbersome way in \cite[Proposition 2.2]{CBPS2018}.

\begin{lemma}\label{rdt}
Denote by $z$ and $z'$ two successive iterates
of the PSD-id (Algorithm \ref{psdid})
where $z\in\spanl\{V\}{\setminus}\{0\}$.
Let the preconditioner $K$
be effectively positive definite. If $z$
is not an eigenvector, then $\rho(z')<\rho(z)$.
\end{lemma}
\begin{proof}
We use the coefficient vectors $c$ and $c'$ defined in Lemma \ref{lmpsdid}.
Since $c'$ is a minimizer of \,$\widetilde{\rho}(\cdot)$ in
$\spanl\{c,\widetilde{K}\widetilde{r}\}$, we get
$\widetilde{\rho}(c')\le\widetilde{\rho}(c)$.
Therein the equality does not hold, since otherwise
$c$ would also be a minimizer of $\widetilde{\rho}(\cdot)$
and thus a Ritz vector in $\spanl\{c,\widetilde{K}\widetilde{r}\}$.
Then the residual $\widetilde{r}=\Lambda c-\widetilde{\rho}(c)c$
would be orthogonal to $\spanl\{c,\widetilde{K}\widetilde{r}\}$
so that $\widetilde{r}^*\widetilde{K}\widetilde{r}=0$.
Subsequently, the positive definiteness of the restricted form
$\widetilde{K}$ of $K$ leads to $\widetilde{r}=0$ and
\[
 0=SV\widetilde{r}=SV\big(\Lambda c-\widetilde{\rho}(c)c\big)
 \stackrel{\eqref{bevp}}{=}HVc-\widetilde{\rho}(c)SVc
 \stackrel{\eqref{rqe}}{=}Hz-\rho(z)Sz,
\]
i.e., $z$ would be an eigenvector. Thus
$\widetilde{\rho}(c')<\widetilde{\rho}(c)$ holds
and implies $\rho(z')<\rho(z)$ by \eqref{rqe}.
\end{proof}

The following lemma
provides a quantitive measure on the quality of 
an effectively positive definite preconditioner.

\begin{lemma}\label{pcd}
Consider an effectively positive definite preconditioner $K$,
its restricted form $\widetilde{K}=V^*SKSV$
and the diagonal matrix $\Lambda_{\nu}=\Lambda-\nu \widetilde{I}$, 
where $\Lambda$ is from \eqref{bevp}, $\widetilde{I}=I_{n-i+1}$,
and $\nu$ is a parameter such that $\nu<\la_i$.
Denote by $\alpha$ and $\beta$ the smallest and 
largest eigenvalues of $\widetilde{K}\Lambda_{\nu}$.
Then $\beta\ge\alpha>0$, and
\begin{equation}\label{pcde}
\big\|\widetilde{I}-\omega\widetilde{K}\Lambda_{\nu}\big\|_{\Lambda_{\nu}}
\le\varepsilon<1,
\end{equation}
where $\omega=2/(\beta+\alpha)$
and $\varepsilon=(\beta-\alpha)/(\beta+\alpha)$.
\end{lemma}
\begin{proof}
The matrices $\widetilde{K}$ and $\Lambda_{\nu}$ are
evidently Hermitian positive definite so that
their square root matrices are available.
By using $\Lambda_{\nu}^{1/2}$, the matrix
$\widetilde{K}\Lambda_{\nu}$ is similar to
$\widehat{K}=\Lambda_{\nu}^{1/2}\widetilde{K}\Lambda_{\nu}^{1/2}$
which is Hermitian positive definite. This shows the positiveness
of all eigenvalues of $\widetilde{K}\Lambda_{\nu}$ and $\widehat{K}$.
Moreover, the norm
$\big\|\widetilde{I}-\omega\widetilde{K}\Lambda_{\nu}\big\|_{\Lambda_{\nu}}
 =\big\|\widetilde{I}-\omega\widehat{K}\big\|_2$
is actually the maximum of $|1-\omega\lambda|$ among all eigenvalues $\lambda$
of $\widehat{K}$. Then the quality condition \eqref{pcde} is verified by
\[0<\alpha\le\lambda\le\beta
 \quad\Rightarrow\quad
 \frac{2\alpha}{\beta+\alpha}\le\omega\lambda\le\frac{2\beta}{\beta+\alpha}
 \quad\Rightarrow\quad
 \frac{\alpha-\beta}{\beta+\alpha}\le1-\omega\lambda
 \le\frac{\beta-\alpha}{\beta+\alpha}.\]
\end{proof}

The following lemma interprets the coefficient vectors from Lemma~\ref{lmpsdid}
as iterates of a PSD iteration for a shifted matrix.

\begin{lemma}\label{lmpsdid1}
With  the diagonal matrix $\Lambda_{\nu}$ from Lemma~\ref{pcd}
and the corresponding Rayleigh quotient
\begin{equation}\label{lrqs}
\widetilde{\rho}_{\nu}: \CC^{n-i+1}{\setminus}\{0\}\to\RR,\quad
\widetilde{\rho}_{\nu}(w)=\frac{w^*\Lambda_{\nu}w}{w^*w},
\end{equation}
the coefficient vector $c'$ of $z'$ for the PSD-id (Algorithm \ref{psdid})
is a minimizer of \,$\widetilde{\rho}_{\nu}(\cdot)$ in
$\spanl\{c,\widetilde{K}\widetilde{r}_{\nu}\}$, 
where $\widetilde{K}=V^*SKSV$ and 
$\widetilde{r}_{\nu}=\Lambda_{\nu} c-\widetilde{\rho}_{\nu}(c)c$.
\end{lemma}
\begin{proof}
For an arbitrary $\widetilde{c}\in\CC^{n-i+1}{\setminus}\{0\}$, it holds that
\[
 \widetilde{\rho}_{\nu}(\widetilde{c})
 =\frac{\widetilde{c}^*\Lambda_{\nu}\widetilde{c}}{\widetilde{c}^*\widetilde{c}}
 =\frac{\widetilde{c}^*\Lambda \widetilde{c}-\nu\widetilde{c}^*\widetilde{c}}
  {\widetilde{c}^*\widetilde{c}}
 =\widetilde{\rho}(\widetilde{c})-\nu.
\]
Thus minimizing $\widetilde{\rho}_{\nu}(\cdot)$ is equivalent to
minimizing $\widetilde{\rho}(\cdot)$. Moreover, the relation
\[
 \widetilde{r}_{\nu}=\Lambda_{\nu} c-\widetilde{\rho}_{\nu}(c)c
 =\big(\Lambda-\nu \widetilde{I}\big)c-\big(\widetilde{\rho}(c)-\nu\big)c
 =\Lambda c-\widetilde{\rho}(c)c=\widetilde{r}
\]
implies $\spanl\{c,\widetilde{K}\widetilde{r}_{\nu}\}
=\spanl\{c,\widetilde{K}\widetilde{r}\}$
so that the statement for $c'$ from Lemma \ref{lmpsdid}
is directly reformulated in terms of
$\widetilde{\rho}_{\nu}(\cdot)$ and $\widetilde{r}_{\nu}$.
\end{proof}

By Lemmas~\ref{pcd} and \ref{lmpsdid1},
the following theorem shows that by a reverse transformation,
the PSD estimate \eqref{estpsd} in Theorem~\ref{thmbpsd} leads to
a sharp single-step estimate of the PSD-id 
based on relations of Rayleigh quotients in \eqref{rqe}.

\begin{theorem}\label{psdidth}
Denote by $z$ and $z'$ two successive iterates
of the PSD-id (Algorithm \ref{psdid})
where $z\in\spanl\{V\}{\setminus}\{0\}$.
Let the preconditioner $K$ be effectively positive definite
with the quality parameter $\varepsilon$ defined in \eqref{pcde}
for any $\nu<\la_i$.

If $\rho(z)\in(\la_j,\la_{j+1})$ for certain $j \ge i$, then
\begin{equation}\label{eq:psdidbd}
 \frac{\rho(z')-\la_j}{\la_{j+1}-\rho(z')}
 \le\left(\frac{\kappa+\varepsilon(2-\kappa)}{(2-\kappa)+\varepsilon\kappa}\right)^2
 \frac{\rho(z)-\la_j}{\la_{j+1}-\rho(z)}
\end{equation}
with
\[\kappa=\left(\frac{\la_j-\nu}{\la_{j+1}-\nu}\right)
 \left(\frac{\la_n-\la_{j+1}}{\la_n-\la_j}\right).\]
The equality in \eqref{eq:psdidbd} is attainable
in the limit case $\rho(z)\to\la_j$
in an invariant subspace associated with
the eigenvalues $\la_j$, $\la_{j+1}$ and $\la_n$.
\end{theorem}
\begin{proof}
We use coefficient vectors $c$ and $c'$ introduced in Lemma \ref{lmpsdid}.
According to Lemma \ref{lmpsdid1},
$c'$ is a minimizer of \,$\widetilde{\rho}_{\nu}(\cdot)$ in
$\spanl\{c,\widetilde{K}\widetilde{r}_{\nu}\}$
concerning the matrix $\Lambda_{\nu}$.

By applying Theorem \ref{thmbpsd} to
\[
 A\to\Lambda_{\nu}, \quad
 x \to c, \quad
 T\to\omega\widetilde{K}, \quad
 \eta\to\widetilde{\rho}_{\nu}(c), \quad
 \eta'\to\widetilde{\rho}_{\nu}(c')
\]
and substituting the eigenvalues, the estimate \eqref{estpsd} is specified as
\[
 \frac{\widetilde{\rho}_{\nu}(c')-(\la_j-\nu)}
  {\,(\la_{j+1}-\nu)-\widetilde{\rho}_{\nu}(c')\,}
 \le\left(\frac{\kappa+\varepsilon(2-\kappa)}
  {(2-\kappa)+\varepsilon\kappa}\right)^2
 \frac{\widetilde{\rho}_{\nu}(c)-(\la_j-\nu)}
  {\,(\la_{j+1}-\nu)-\widetilde{\rho}_{\nu}(c)\,}
\]
with
\[
 \kappa=\frac{(\la_j-\nu)\big((\la_n-\nu)-(\la_{j+1}-\nu)\big)}
  {(\la_{j+1}-\nu)\big((\la_n-\nu)-(\la_j-\nu)\big)}
 =\left(\frac{\la_j-\nu}{\la_{j+1}-\nu}\right)
 \left(\frac{\la_n-\la_{j+1}}{\la_n-\la_j}\right).
\]
Therein the terms $\widetilde{\rho}_{\nu}(c)$ and
$\widetilde{\rho}_{\nu}(c')$ coincide with
$\rho(z)-\nu$ and $\rho(z')-\nu$ due to the relation
\[
 \widetilde{\rho}_{\nu}(\widetilde{c})
 =\widetilde{\rho}(\widetilde{c})-\nu
 \stackrel{\eqref{rqe}}{=}\rho(\widetilde{z})-\nu.
\]
Thus \eqref{eq:psdidbd} is shown. Furthermore,
the sharpness statement in Theorem \ref{thmbpsd}
is specified for the limit case $\widetilde{\rho}_{\nu}(c)\to\la_j-\nu$
and the matrix $\Lambda_{\nu}$. The corresponding invariant subspace
is associated with the eigenvalues $\la_j-\nu$, $\la_{j+1}-\nu$ and $\la_n-\nu$.
Denoting this subspace by $\widetilde{\mathcal{C}}$, then
$\widetilde{\mathcal{Z}}=V\widetilde{\mathcal{C}}$ is an invariant subspace
of $(H,S)$ associated with the eigenvalues $\la_j$, $\la_{j+1}$ and $\la_n$
due to \eqref{bevp}, and the limit case $\widetilde{\rho}_{\nu}(c)\to\la_j-\nu$
is converted into $\rho(z)\to\la_j$.
\end{proof}

\begin{remark}\label{rmrdt}
The assumption $\rho(z)\in(\la_j,\la_{j+1})$ in Theorem~\ref{psdidth}
does not cover the case that $\rho(z)$ is equal to $\la_j$ or $\la_{j+1}$.
Applying Lemma~\ref{rdt} to e.g. $\rho(z)=\la_j$ provides 
the following supplement:
if $z$ is an eigenvector, then the iteration is simply terminated;
otherwise Lemma~\ref{rdt} ensures that
$\rho(z)$ is smaller in the next step
and can match the assumption $\rho(z)\in(\la_j,\la_{j+1})$
for a smaller index $j$ so that Theorem~\ref{psdidth} is applicable.
In summary, $\rho(z)$ either converges to an eigenvalue $\la_j$ with $j>i$
or reaches the interval $(\la_i,\la_{i+1})$ in the final phase. 
In the latter (and usual) case, two possible phenomena 
can be interpreted by the ratio $(\la_i-\nu)/(\la_{i+1}-\nu)$
from the convergence bound:
the convergence rate is deteriorated for $\la_i\approx\la_{i+1}$;
for well-separated $\la_i$ and $\la_{i+1}$,
a fast convergence can be obtained by
(proper approximations of) the shift-invert preconditioner
$K=(H-\sigma S)^{-1}$ with $\sigma\approx\la_i$
since $\kappa\to0$ for $\nu=\sigma\to\la_i$.
In \cite[Section 5]{CBPS2018}, it is shown that an efficient shift 
$\sigma$ can be chosen
from the interval $(\la_{i-1},\la_i)$, e.g., by initially setting $\sigma$ slightly larger than
the computed $\la_{i-1}$ and then enlarging it with a weighted mean
of $\la_{i-1}$ and the current approximation of $\la_i$.
\end{remark}

\begin{remark}\label{comparison}
In comparison to \cite[Theorem 3.2]{CBPS2018}, the current approximate eigenvalue
$\rho(z)$ in Theorem \ref{psdidth} is located in an arbitrary eigenvalue interval
so that the statement is much more flexible.
Moreover, the bound in \eqref{eq:psdidbd}
has a simpler form where only one technical term is used,
namely the quality parameter $\varepsilon$ of preconditioning.
With a dynamic shift $\sigma$ approximating $\la_i$ from below,
the parameter $\kappa$ in Theorem~\ref{psdidth}
with $\nu=\sigma$ can be close to zero
in the final phase and indicates a superlinear convergence.
The limit case $\kappa\to0$ corresponds to an optimal
shift-invert preconditioner which allows a one-step convergence.
In addition, \cite[Theorem 4.1]{CBPS2018} can be improved by
Theorem~\ref{psdidth} with the convergence factor $\kappa^2/(2-\kappa)^2$
arising from the special case $\varepsilon=0$, i.e., $K=(H-\sigma S)^{-1}$.
\end{remark}

\subsection{Sharp single-step estimate of the BPSD-id} \label{sec:bpsdid}

Let us now analyze the evolution of Ritz values of the BPSD-id
within two successive subspace iterates.
We first interpret the coefficient matrices of the BPSD-id from 
Lemma~\ref{lmbpsdid}
as iterates of a BPSD iteration for the shifted matrix $\Lambda_{\nu}$
introduced in Lemma~\ref{pcd}.

\begin{lemma}\label{lmbpsdid1}
Denote by $Z$ and $Z'$ two successive iterates
of the BPSD-id (Algorithm \ref{bpsdid})
where $Z$ has full rank and all its $\widetilde{k}$ columns
belong to $\spanl\{V\}$. With the diagonal matrix $\Lambda_{\nu}$
from Lemma \ref{pcd}, the columns
of the coefficient matrix $C=V^*SZ$ are orthonormal Ritz vectors
of $\Lambda_{\nu}$ in $\spanl\{C\}$.
Moreover, by using the corresponding Ritz value matrix 
$\Theta_{\nu}=C^*\Lambda_{\nu}C$, the columns 
of the coefficient matrix $C'=V^*SZ'$
are orthonormal Ritz vectors associated with
the $\widetilde{k}$ smallest Ritz values of $\Lambda_{\nu}$ in
$\spanl\{C,\widetilde{K}\widetilde{R}_{\nu}\}$
for $\widetilde{K}=V^*SKSV$ and 
$\widetilde{R}_{\nu}=\Lambda_{\nu}C - C\Theta_{\nu}$.
\end{lemma}
\begin{proof}
The statements are verified by using Lemma \ref{lmbpsdid}
and the transformations
\[\Theta_{\nu}
 =C^*\Lambda_{\nu}C
 =C^*\big(\Lambda-\nu\widetilde{I}\,\big)C
 =\Theta-\nu I_{\widetilde{k}},
\]
and the fact that 
\[
 \widetilde{R}_{\nu}
 =\big(\Lambda-\nu \widetilde{I}\,\big)C-
C\big(\Theta-\nu I_{\widetilde{k}}\big)
 =\Lambda C-C\Theta=\widetilde{R}
\]
implies that 
\[
 \spanl\{C,\widetilde{K}\widetilde{R}_{\nu}\}
 =\spanl\{C,\widetilde{K}\widetilde{R}\}.
\] 
\end{proof}

The following lemma shows a strict reduction of Ritz values
concerning the coefficient subspaces $\spanl\{C\}$ and $\spanl\{C'\}$.

\begin{lemma}\label{brdt}
Let the preconditioner $K$ of the BPSD-id
be effectively positive definite.
Following Lemma~\ref{lmbpsdid1},
denote by $\varphi_1\le\cdots\le\varphi_{\widetilde{k}}$
and $\varphi'_1\le\cdots\le\varphi'_{\widetilde{k}}$
the Ritz values of $\Lambda_{\nu}$
in $\spanl\{C\}$ and $\spanl\{C'\}$, respectively.
If $\spanl\{C\}$ contains no eigenvectors of $\Lambda_{\nu}$,
then $\varphi'_t<\varphi_t$ holds for each $t\in\{1,\ldots,\widetilde{k}\}$.
\end{lemma}
\begin{proof}
Let us apply Lemma~\ref{pcd} to $K$, and define the auxiliary matrix
\[\widetilde{C}=C-\omega\widetilde{K}\widetilde{R}_{\nu}.\]
Then $\widetilde{C}$ has full rank since otherwise there would exist
a $g\in\CC^{\widetilde{k}}{\setminus}\{0\}$ satisfying
$0=\widetilde{C}g
=Cg-\omega\widetilde{K}(\Lambda_{\nu}C - C\Theta_{\nu})g$, i.e.,
\[
 \Lambda_{\nu}^{-1}C\Theta_{\nu}g
 =(\widetilde{I}-\omega\widetilde{K}\Lambda_{\nu})
 (\Lambda_{\nu}^{-1}C\Theta_{\nu}g-Cg).
\]
Moreover, the condition $\nu<\la_i$ from Lemma \ref{pcd}
ensures that $\Lambda_{\nu}^{-1}C\Theta_{\nu}$ has full rank
so that $\Lambda_{\nu}^{-1}C\Theta_{\nu}g\neq0$.
Applying \eqref{pcde} implies
\begin{equation}\label{ctd}
 \|\Lambda_{\nu}^{-1}C\Theta_{\nu}g\|_{\Lambda_{\nu}}
 \le\|\widetilde{I}-\omega\widetilde{K}\Lambda_{\nu}\|_{\Lambda_{\nu}}
 \|\Lambda_{\nu}^{-1}C\Theta_{\nu}g-Cg\|_{\Lambda_{\nu}}
 <\|\Lambda_{\nu}^{-1}C\Theta_{\nu}g-Cg\|_{\Lambda_{\nu}}, 
\end{equation}
where $\|\Lambda_{\nu}^{-1}C\Theta_{\nu}g-Cg\|_{\Lambda_{\nu}}=0$
is excluded due to the first inequality and
$\Lambda_{\nu}^{-1}C\Theta_{\nu}g\neq0$.
However, the orthogonality
\[
 (Cg)^*\Lambda_{\nu}(\Lambda_{\nu}^{-1}C\Theta_{\nu}g-Cg)
 =g^*C^*C\Theta_{\nu}g-g^*C^*\Lambda_{\nu}Cg=0
\]
leads to
\[
 \|\Lambda_{\nu}^{-1}C\Theta_{\nu}g\|_{\Lambda_{\nu}}^2
 =\|Cg\|_{\Lambda_{\nu}}^2
 +\|\Lambda_{\nu}^{-1}C\Theta_{\nu}g-Cg\|_{\Lambda_{\nu}}^2
 \ge\|\Lambda_{\nu}^{-1}C\Theta_{\nu}g-Cg\|_{\Lambda_{\nu}}^2
\]
which contradicts \eqref{ctd}. Thus $\widetilde{C}$ has full rank.

Consequently, there are $\widetilde{k}$ Ritz values in $\spanl\{\widetilde{C}\}$.
We denote them by $\widetilde{\varphi}_1\le\cdots\le\widetilde{\varphi}_{\widetilde{k}}$.
Then $\varphi'_t\le\widetilde{\varphi}_t$ holds for each $t\in\{1,\ldots,{\widetilde{k}}\}$
due to $\spanl\{\widetilde{C}\}\subseteq\spanl\{C,\widetilde{K}\widetilde{R}_{\nu}\}$
and the Courant-Fischer principles. For proving $\varphi'_t<\varphi_t$,
it remains to be shown 
\begin{equation} \label{eq:phibd} 
\widetilde{\varphi}_t<\varphi_t.
\end{equation}
For this purpose, we use a submatrix $E_t$ of $I_{\widetilde{k}}$ such that
the columns of $CE_t$ are Ritz vectors associated with
the Ritz values $\varphi_1\le\cdots\le\varphi_t$.
Then $\varphi_t$ is the largest Ritz value in $\spanl\{CE_t\}$.
Correspondingly, we consider the largest Ritz value
$\widehat{\varphi}_t$ in $\spanl\{\widetilde{C}E_t\}$.
The relation $\spanl\{\widetilde{C}E_t\}\subseteq\spanl\{\widetilde{C}\}$
ensures $\widetilde{\varphi}_t\le\widehat{\varphi}_t$.
In addition, we use a Ritz vector $\widehat{c}=\widetilde{C}E_t\widehat{g}$
associated with $\widehat{\varphi}_t$ and the auxiliary vector
\[
 c=C\Theta_{\nu}E_t\widehat{g}
 =CE_t\Theta_{\nu,t}\,\widehat{g}
 \quad\mbox{with}\quad
 \Theta_{\nu,t}=\diag(\varphi_1,\ldots,\varphi_t).
\]
Then by using the definitions of $\widetilde{C}$ and $\widetilde{R}_{\nu}$, 
we have
\[
 \Lambda_{\nu}^{-1}c-\widehat{c}
 =\Lambda_{\nu}^{-1}C\Theta_{\nu}E_t\widehat{g}
 -\big(C-\omega\widetilde{K}(\Lambda_{\nu}C - C\Theta_{\nu})\big)E_t\widehat{g}
 =(\widetilde{I}-\omega\widetilde{K}\Lambda_{\nu})
 (\Lambda_{\nu}^{-1}c-CE_t\widehat{g})
\]
so that
\[
 \|\Lambda_{\nu}^{-1}c-\widehat{c}\|_{\Lambda_{\nu}}
 \stackrel{\eqref{pcde}}{\le}
 \varepsilon\|\Lambda_{\nu}^{-1}c-CE_t\widehat{g}\|_{\Lambda_{\nu}}.
\]
Therein $\|\Lambda_{\nu}^{-1}c-CE_t\widehat{g}\|_{\Lambda_{\nu}}$
further fulfills
\[
 \|\Lambda_{\nu}^{-1}c-CE_t\widehat{g}\|_{\Lambda_{\nu}}
 \le\|\Lambda_{\nu}^{-1}c-\varphi^{-1}c\|_{\Lambda_{\nu}}
\]
for $\varphi=\widetilde{\rho}_{\nu}(c)$ with \eqref{lrqs} since the orthogonality
\[
 (CE_t\widetilde{g})^*\Lambda_{\nu}
 (\Lambda_{\nu}^{-1}c-CE_t\widehat{g})
 =\widetilde{g}^*E_t^*C^*C\Theta_{\nu}E_t\widehat{g}
 -\widetilde{g}^*E_t^*C^*\Lambda_{\nu}CE_t\widehat{g}=0
\]
for $\widetilde{g}=\widehat{g}-\varphi^{-1}\Theta_{\nu,t}\widehat{g}$ leads to
\[
 \begin{split}
 \|\Lambda_{\nu}^{-1}c-CE_t\widehat{g}\|_{\Lambda_{\nu}}
 &\le\big(\|\Lambda_{\nu}^{-1}c-CE_t\widehat{g}\|_{\Lambda_{\nu}}^2
 +\|CE_t\widetilde{g}\|_{\Lambda_{\nu}}^2\big)^{1/2} \\[1ex]
 &=\|\Lambda_{\nu}^{-1}c
 -CE_t(\varphi^{-1}\Theta_{\nu,t}\widehat{g})\|_{\Lambda_{\nu}}
 =\|\Lambda_{\nu}^{-1}c-\varphi^{-1}c\|_{\Lambda_{\nu}}.
 \end{split}
\]
If $\spanl\{C\}$ contains no eigenvectors of $\Lambda_{\nu}$, then
the difference $\Lambda_{\nu}^{-1}c-\varphi^{-1}c$ is nonzero since otherwise
$c$ would be an eigenvector of $\Lambda_{\nu}$. Consequently, it holds that
\[
 \|\Lambda_{\nu}^{-1}c-\widehat{c}\|_{\Lambda_{\nu}}^2
 <\|\Lambda_{\nu}^{-1}c-\varphi^{-1}c\|_{\Lambda_{\nu}}^2.
\]
By the definition of the norm $\|\cdot\|_{\Lambda_{\nu}}$, 
we have 
\[ 
 \widehat{c}^*\Lambda_{\nu}\widehat{c}
 -c^*\widehat{c}-\widehat{c}^*c
 <\varphi^{-2}c^*\Lambda_{\nu}c-2\varphi^{-1}c^*c
 =-\varphi^{-1}c^*c. 
\]
Multiplying both sides by $\varphi^{-1}$ yields
\[ 
 \varphi^{-1}\widehat{c}^*\Lambda_{\nu}\widehat{c}
 <\varphi^{-1}(c^*\widehat{c}+\widehat{c}^*c)-\varphi^{-2}c^*c
 =\widehat{c}^*\widehat{c}-\|\widehat{c}-\varphi^{-1}c\|_2^2
 \le\widehat{c}^*\widehat{c}
\]
so that $\widetilde{\rho}_{\nu}(\widehat{c})
=(\widehat{c}^*\Lambda_{\nu}\widehat{c})/(\widehat{c}^*\widehat{c})
<\varphi$. Thus the inequality~\eqref{eq:phibd} is verified by
\[ 
 \widetilde{\varphi}_t\le\widehat{\varphi}_t
 =\widetilde{\rho}_{\nu}(\widehat{c})
 <\varphi=\widetilde{\rho}_{\nu}(c)
 =\widetilde{\rho}_{\nu}(CE_t\Theta_{\nu,t}\,\widehat{g})
 \le\varphi_t.
\]
\end{proof}

Lemma \ref{brdt} can easily be adapted to successive iterates of the BPSD-id.
Therein the relation \eqref{rve} applied to $C$, $Z$ and $C'$, $Z'$
(introduced in Lemma \ref{lmbpsdid1})
shows that $\varphi_t+\nu$ and $\varphi'_t+\nu$
are Ritz values of $(H,S)$ in $\spanl\{Z\}$ and $\spanl\{Z'\}$, respectively.
Therefore a strict reduction of Ritz values occurs
if $\spanl\{Z\}$ contains no eigenvectors of $(H,S)$.
This actually generalizes Lemma~\ref{rdt} to the BPSD-id
where the proof is however not a direct generalization
of that of  Lemma~\ref{rdt}.

The following theorem includes
sharp single-step estimates on the convergence of the BPSD-id.

\begin{theorem}\label{bpsdidth}
Denote by $Z$ and $Z'$ two successive iterates
of the BPSD-id (Algorithm \ref{bpsdid})
where $Z$ has full rank and all its $\widetilde{k}$ columns
belong to $\spanl\{V\}$. 
Let the preconditioner $K$ be effectively positive definite
with the quality parameter $\varepsilon$ defined in \eqref{pcde}.
Consider the Ritz values
$\theta_1\le\cdots\le\theta_{\widetilde{k}}$
and $\theta'_1\le\cdots\le\theta'_{\widetilde{k}}$
in $\spanl\{Z\}$ and $\spanl\{Z'\}$, respectively.

If $\theta_t\in(\la_j,\la_{j+1})$
for $t\in\{1,\ldots,\widetilde{k}\}$ and certain $j \ge i{-}1{+}t$, then
\begin{equation}\label{eq:bpsdidbd}
 \frac{\theta'_t-\la_j}{\la_{j+1}-\theta'_t}
 \le\left(\frac{\kappa+\varepsilon(2-\kappa)}{(2-\kappa)+\varepsilon\kappa}\right)^2
 \frac{\theta_t-\la_j}{\la_{j+1}-\theta_t}
\end{equation}
with
\[\kappa=\left(\frac{\la_j-\nu}{\la_{j+1}-\nu}\right)
 \left(\frac{\la_n-\la_{j+1}}{\la_n-\la_j}\right).\]
The equality in \eqref{eq:bpsdidbd} is attainable
in the limit case $\theta_t\to\la_j$
in an invariant subspace associated with
the eigenvalues $\la_j$, $\la_{j+1}$ and $\la_n$.
\end{theorem}
\begin{proof}
We use the coefficient matrices $C$ and $C'$ 
of $Z$ and $Z'$ introduced in Lemma \ref{lmbpsdid}.
Following the relation \eqref{rve},
Lemma \ref{lmbpsdid1} and Lemma \ref{brdt},
we apply Theorem \ref{thmbpsd} to
\[
 A \to\Lambda_{\nu}, \quad
 X \to C, \quad
 T\to\omega\widetilde{K}, \quad
 \eta_t\to\varphi_t=\theta_t-\nu, \quad
 \eta'_t\to\varphi'_t=\theta'_t-\nu.
\]
This results in
\[
 \frac{(\theta'_t-\nu)-(\la_j-\nu)}
  {\,(\la_{j+1}-\nu)-(\theta'_t-\nu)\,}
 \le\left(\frac{\kappa+\varepsilon(2-\kappa)}
  {(2-\kappa)+\varepsilon\kappa}\right)^2
 \frac{(\theta_t-\nu)-(\la_j-\nu)}
  {\,(\la_{j+1}-\nu)-(\theta_t-\nu)\,}
\]
with
\[
 \kappa=\frac{(\la_j-\nu)\big((\la_n-\nu)-(\la_{j+1}-\nu)\big)}
  {(\la_{j+1}-\nu)\big((\la_n-\nu)-(\la_j-\nu)\big)}
 =\left(\frac{\la_j-\nu}{\la_{j+1}-\nu}\right)
 \left(\frac{\la_n-\la_{j+1}}{\la_n-\la_j}\right)
\]
and further yields the estimate \eqref{eq:bpsdidbd}
including the sharpness statement.
\end{proof}

Similarly to Remark~\ref{rmrdt},
we can additionally discuss the case $\theta_t=\la_j$
based on Lemma~\ref{brdt}. This yields
\[
 \theta'_t-\nu=\varphi'_t<\varphi_t=\theta_t-\nu
 \quad\Rightarrow\quad
 \theta'_t<\theta_t
\]
so that the limit of the Ritz value $\theta_t$ is either $\la_{i-1+t}$ or 
some $\la_j$ with $j > i{-}1{+}t$.

We note that the convergence factor in \eqref{eq:bpsdidbd} 
is only meaningful for
well-separated eigenvalues and cannot predict the so-called cluster robustness
of block iterations. Section~\ref{sec:multi} serves to fill this theoretical gap.

\subsection{Larger shifts}  \label{sec:sigma}

The estimates in Theorems~\ref{psdidth} and~\ref{bpsdidth}
are concerned with effectively positive definite preconditioners
such as shift-invert preconditioners of the form $K=(H-\sigma S)^{-1}$
as well as their approximations. Setting $\sigma<\la_i$
easily ensures the effectively positive definiteness of $K$; see \eqref{irf}.
However, the practical implementation of the PSD-id in \cite[Section 5]{CBPS2018}
uses the shift $\sigma=\rho(z)$ after $\rho(z)$ is sufficiently close to $\la_i$,
i.e., $\sigma$ is larger than $\la_i$. Nevertheless,
the convergence analysis therein is still based on some estimates
which actually treat the case that $\sigma$ is slightly smaller than $\la_i$.
Thus the resulting bounds can only be applied in an asymptotic way.
This section presents a direct analysis together with an extension to the BPSD-id.

We begin with the condition $\la_i<\sigma<(\la_i+\la_{i+1})/2$ in the PSD-id 
and the exact shift-invert preconditioner $K=(H-\sigma S)^{-1}$.

\begin{theorem}\label{psdidrq}
Let $z$ be the current iterate of the PSD-id (Algorithm \ref{psdid}),
and $z'$ the next iterate. Assume that the columns $u_1,\ldots,u_{i-1}$ of $U$
are $S$-orthonormal eigenvectors associated with the $i-1$ smallest eigenvalues
$\lambda_1\le\cdots\le\lambda_{i-1}$, define $V=[u_i,\ldots,u_n]$
such that $[U,V]=[u_1,\ldots,u_n]$ is an $S$-orthonormal eigenbasis
associated with $\lambda_1\le\cdots\le\lambda_n$, and consider the case
$K=(H-\sigma S)^{-1}$ 
for a shift $\sigma\in\big(\la_i,\,(\la_i+\la_{i+1})/2\big)$.

If $\rho(z)\neq\sigma$ and $\rho(z)\in(\la_i,\la_{i+1})$, then
\begin{equation}\label{psdidrqe}
 \frac{\rho(z')-\la_i}{\la_{i+1}-\rho(z')}
 \le\left(\frac{\kappa}{2-\kappa}\right)^2
 \frac{\rho(z)-\la_i}{\la_{i+1}-\rho(z)}
\end{equation}
with
\[\kappa=\left(\frac{\la_i-\sigma}{\la_{i+1}-\sigma}\right)
 \left(\frac{\la_n-\la_{i+1}}{\la_n-\la_i}\right).\]
\end{theorem}
\begin{proof}
The condition $\rho(z)\neq\sigma$ excludes a stagnation caused by
\[p=\big(H-\rho(z)S\big)^{-1}(Hz-\rho(z)Sz)=z
 \quad\Rightarrow\quad \spanl\{U,z,p\}=\spanl\{U,z\}.\]
For proving \eqref{psdidrqe}, we use a restricted formulation of
the considered iteration with respect to $z=Vc$ and $z'=Vc'$, namely,
\begin{equation}\label{psdida}
c'=\mathrm{RR}
\big(\spanl\{c,\,(\Lambda - \sigma \widetilde{I})^{-1}c\}\big)
\end{equation}
where the Rayleigh-Ritz procedure $\mathrm{RR}(\cdot)$
extracts an orthonormal Ritz vector of $\Lambda$ associated with the
smallest Ritz value; cf.~Lemma \ref{lmpsdid}.
Evidently, \eqref{psdida} is an acceleration of
\[\widehat{c}=(\Lambda - \sigma \widetilde{I})^{-1}c+\beta c\]
for arbitrary $\beta\in\RR$. By using the auxiliary matrix $A=-\Lambda$
together with its eigenvalues $\alpha_1\ge\cdots\ge\alpha_{n-i+1}$
and the corresponding Rayleigh quotient $\alpha(\cdot)$, we get
\[\alpha_j=-\la_{i-1+j},\quad
 \alpha(c)=-\widetilde{\rho}(c)=-\rho(z)>-\la_{i+1}=\alpha_2,\quad
 \alpha(\widehat{c})=-\widetilde{\rho}(\widehat{c}),\]
\[\mbox{and}\quad
 \widehat{c}=f(A)c
 \quad\mbox{with}\quad
 f(\eta)=(-\eta-\sigma)^{-1}+\beta.\]
Subsequently, we choose
$\beta=-\tfrac12\big((\la_{i+1}-\sigma)^{-1}+(\la_n-\sigma)^{-1}\big)$
so that
\[|f(\alpha_1)|>|f(\alpha_2)|\ge|f(\alpha_j)|
 \quad\forall\ j\in\{2,\,\ldots,\,n{\,-\,}i{\,+\,}1\}\]
holds (by elementary comparison). Then $\widehat{c}=f(A)c$
can be analyzed as the power method for a matrix function
\cite[Section 1.1]{KNY1987}. In particular,
the estimate \cite[(1.9)]{KNY1987} is applicable due to $\alpha(c)>\alpha_2$,
and implies
\[\frac{\alpha_1-\alpha(\widehat{c})}{\alpha(\widehat{c})-\alpha_2}
 \le\left(\frac{|f(\alpha_2)|}{|f(\alpha_1)|}\right)^2\
 \frac{\alpha_1-\alpha(c)}{\alpha(c)-\alpha_2}\]
which is equivalent to
\[\frac{\widetilde{\rho}(\widehat{c})-\la_i}
 {\la_{i+1}-\widetilde{\rho}(\widehat{c})}
 \le\left(\frac{(\la_{i+1} - \sigma)^{-1}+\beta}{(\la_i - \sigma)^{-1}+\beta}\right)^2
 \ \frac{\widetilde{\rho}(c)-\la_i}{\la_{i+1}-\widetilde{\rho}(c)}.\]
This leads to \eqref{psdidrqe} by inserting $\beta$ and using
$\widetilde{\rho}(\widehat{c})\ge\widetilde{\rho}(c')=\rho(z')$, \
$\widetilde{\rho}(c)=\rho(z)$.
\end{proof}

Theorem~\ref{psdidrq} relaxes the condition $\sigma<\la_i$.
A reasonable interval for selecting the shift
is $\big(\la_{i-1},\,(\la_i+\la_{i+1})/2\big)$ where the lower bound
can be determined by the previous outer step and the upper bound
can be detected by the Rayleigh quotient and the residual;
cf.~\cite[Section 5.2]{NTY2002}.
The excluded case $\rho(z)=\sigma$ in Theorem~\ref{psdidrq}
can be treated by solving the restricted linear system
\begin{equation}\label{jde}
\widetilde{Q}^*(H-\sigma S)\widetilde{Q}p=r, \quad p\in\spanl\{U,z\}^{\perp_S}
\end{equation}
with the $S$-orthogonal projector $\widetilde{Q}$ onto $\spanl\{U,z\}^{\perp_S}$.
Then the subspace $\spanl\{U,z,p\}$
can be shown to contain the vector $(H-\sigma S)^{-1}Sz$ as in \cite[Section 4]{NTY2003}
so that the restricted formulation \eqref{psdida} is applicable.
Therefore the estimate \eqref{psdidrqe} with $j=i$ additionally holds 
for $\rho(z)=\sigma$. Moreover, the special $\beta$ used for deriving \eqref{psdidrqe}
is determined by minimizing $|f(\alpha_2)|/|f(\alpha_1)|$ among $\beta\in\RR$.
Thus setting $\beta=0$ yields a less accurate estimate than \eqref{psdidrqe}, namely,
\begin{equation}\label{cubic}
 \frac{\rho(z')-\la_i}{\la_{i+1}-\rho(z')}
 \le\left(\frac{\rho(z)-\la_i}{\la_{i+1}-\rho(z)}\right)^3.
\end{equation}
In this sense, \eqref{psdidrqe} indicates a ``supercubic'' convergence.

For the PSD-id with inexact shift-invert preconditioners which are not
effectively positive definite, an alternative of the quality parameter
$\varepsilon$ from \eqref{pcde} is required so that \eqref{psdidrqe}
or \eqref{cubic} can be generalized as an estimate like \eqref{eq:psdidbd}.
In particular, for generalizing \eqref{cubic},
an appropriate parameter can be constructed
with respect to the restricted linear system \eqref{jde}.
Therein the restriction of $H-\sigma S$ to
$\spanl\{U,z\}^{\perp_S}$ is positive definite;
cf.~\cite[Lemma 3.1]{NTY2003} or the following simpler verification:
For an arbitrary $w\in\spanl\{U,z\}^{\perp_S}$ with $\|w\|_{S}=1$,
the matrix $W=[U,z,w]$ fulfills $W^*SW=I_{i+1}$ so that
the Ritz values of $(H,S)$ in $\spanl\{W\}$ are eigenvalues of
$W^*HW$. The trace of $W^*HW$ reads
\[\tau=\rho(u_1)+\cdots+\rho(u_{i-1})+\rho(z)+\rho(w)=\la_1+\cdots+\la_{i-1}+\rho(z)+\rho(w).\]
Moreover, the Courant-Fischer principles ensure that $\tau$ is at least
$\la_1+\cdots+\la_{i+1}$. Thus
\[\begin{split}
 &\rho(z)+\rho(w)\ge\la_i+\la_{i+1}\quad\Rightarrow\quad
 \rho(w)\ge\la_i+\la_{i+1}-\sigma\\[1ex]
 &\Rightarrow\quad
 w^*(H-\sigma S)w=\rho(w)-\sigma\ge\la_i+\la_{i+1}-2\sigma>0.
\end{split}\]
Correspondingly, the restricted formulation of PSD-id within $\spanl\{V\}$
includes the restriction of $\Lambda - \sigma \widetilde{I}$ to $\spanl\{c\}^{\perp}$
which is also positive definite. Then the accuracy of the approximate solution $p$
of \eqref{jde} or its coefficient vector $V^*Sp$ can be interpreted by a vector norm
induced by the restriction of $H-\sigma S$ or $\Lambda - \sigma \widetilde{I}$;
cf.~\cite[Theorem 3.2]{NTY2003}. A sufficiently accurate $p$ is characterized 
by a quality parameter $\varepsilon\in[0,1)$ so that the estimate \eqref{cubic}
can be generalized as
\begin{equation}\label{cubic1}
 \frac{\rho(z')-\la_i}{\la_{i+1}-\rho(z')}
 \le\left(\frac{\rho(z)-\la_i+\varepsilon\big(\la_{i+1}-\rho(z)\big)}
 {\la_{i+1}-\rho(z)+\varepsilon\big(\rho(z)-\la_i\big)}\right)^2
 \frac{\rho(z)-\la_i}{\la_{i+1}-\rho(z)}.
\end{equation}

Now we turn to the BPSD-id and analyze its convergence
for the exact shift-inverse preconditioning. 

\begin{theorem}\label{bpsdidrq}
Let $Z$ be the current iterate of the BPSD-id (Algorithm \ref{bpsdid}),
and $Z'$ the next iterate. Assume that the columns $u_1,\ldots,u_{i-1}$ of $U$
are $S$-orthonormal eigenvectors associated with the $i-1$ smallest eigenvalues
$\lambda_1\le\cdots\le\lambda_{i-1}$, define $V=[u_i,\ldots,u_n]$
such that $[U,V]=[u_1,\ldots,u_n]$ is an $S$-orthonormal eigenbasis
associated with $\lambda_1\le\cdots\le\lambda_n$, and consider the case
$K=(H-\sigma S)^{-1}$ for certain $\sigma\in\big(\la_i,\,(\la_i+\la_{i+1})/2\big)$.
Denote by $\theta_1\le\cdots\le\theta_{\widetilde{k}}$
and $\theta'_1\le\cdots\le\theta'_{\widetilde{k}}$
the Ritz values in $\spanl\{Z\}$ and $\spanl\{Z'\}$, respectively.

If $\theta_1\neq\sigma$ and $\theta_1\in(\la_i,\la_{i+1})$, then
\begin{equation}\label{bpsdidrqe}
 \frac{\theta'_1-\la_i}{\la_{i+1}-\theta'_1}
 \le\left(\frac{\kappa}{2-\kappa}\right)^2
 \frac{\theta_1-\la_i}{\la_{i+1}-\theta_1}
\end{equation}
with
\[\kappa=\left(\frac{\la_i-\sigma}{\la_{i+1}-\sigma}\right)
 \left(\frac{\la_n-\la_{i+1}}{\la_n-\la_i}\right).\]
 
If $\theta_t\in(\la_j,\la_{j+1})$
for $t\in\{2,\ldots,\widetilde{k}\}$ and $j = i{\,-\,}1{\,+\,}t$, then
\begin{equation}\label{bpsdidrqe1}
 \frac{\theta'_t-\la_j}{\la_{j+1}-\theta'_t}
 \le\left(\frac{\kappa}{2-\kappa}\right)^2
 \frac{\theta_t-\la_j}{\la_{j+1}-\theta_t}
\end{equation}
with
\[\kappa=\left(\frac{\la_j-\sigma}{\la_{j+1}-\sigma}\right)
 \left(\frac{\la_n-\la_{j+1}}{\la_n-\la_j}\right).\]
\end{theorem}
\begin{proof}
By applying the PSD-id to be a Ritz vector $z_1$ in $\spanl\{Z\}$
associated with $\theta_1$, the next iterate $z'_1$ is contained in $\spanl\{Z'\}$.
Then the estimate \eqref{psdidrqe} for $z=z_1$ and $z'=z'_1$ leads to
\eqref{bpsdidrqe} due to $\rho(z'_1)\ge\theta'_1$ and $\rho(z_1)=\theta_1$.

The derivation of \eqref{bpsdidrqe1} for $t=\widetilde{k}$ is analogous to the proof
of Theorem \ref{psdidrq} by using the generalization
\cite[(2.22)]{KNY1987} of \cite[(1.9)]{KNY1987} 
concerning the block form of an abstract power method.
Moreover, for $t\in\{2,\ldots,\widetilde{k}{\,-\,}1\}$, we can adapt this derivation
to the BPSD-id applied to the subset $\mathcal{Z}_t=\spanl\{z_1,\ldots,z_t\}$
of $\spanl\{Z\}$ spanned by Ritz vectors associated with $\theta_1,\ldots,\theta_t$.
This yields an intermediate estimate for the largest ($t$\,th smallest) Ritz value
$\widetilde{\theta}_t$ in the next iterate $\mathcal{Z}'_t$ of $\mathcal{Z}_t$.
Subsequently, \eqref{bpsdidrqe1} is proved by considering that
$\mathcal{Z}'_t$ is a subset of $\spanl\{Z'\}$ so that $\widetilde{\theta}_t\ge\theta'_t$.
\end{proof}

For the BPSD-id with inexact shift-invert preconditioners,
we particularly consider that
the current $\theta_1$ fulfills $\theta_1<(\la_i+\la_{i+1})/2$
and is used as the next shift $\sigma$. Then the estimates \eqref{cubic}
and \eqref{cubic1} can be adapted to $\theta_1$. It is also remarkable that
by using the positive definiteness of the restriction of $H-\sigma S$,
Theorem \ref{bpsdidth} can be modified
for discussing the deviation between the Ritz values in the subspace iterates
and the Ritz values in the larger subspace $\spanl\{U,Z(:,1)\}^{\perp_S}$.
The latter ones get closer to the eigenvalues $\la_{i+1},\ldots,\la_n$ for
the decreasing $\theta_1$ toward $\la_i$. Then the corresponding estimates turn
into those in Theorem \ref{bpsdidth} with the index update $i \leftarrow i{+1}$.

\section{Multi-step estimates on the cluster robustness} \label{sec:multi} 

A well-known feature of block eigensolvers is the cluster robustness,
i.e., fast convergence toward clustered eigenvalues can be guaranteed
by sufficiently large block sizes.
A classical estimate of the subspace iteration 
has been presented by Parlett \cite[Section 14.4]{PAR1980}.
Therein the block inverse iteration (also called inverse subspace iteration)
$\spanl\{X^{(\ell+1)}\}=\spanl\{A^{-1}X^{(\ell)}\}$ for
a real symmetric positive definite matrix $A$ with eigenvalues
$\alpha_1\le\cdots\le\alpha_n$ is investigated. The convergence is measured
by the angle between an eigenvector and the current subspace.
The resulting bound contains the term $(\alpha_j/\alpha_{m+1})^k$
for $k$ steps with the block size $m$.
A Ritz value estimate with the same ratio
$\alpha_j/\alpha_{m+1}$ can be derived based on the analysis of
the block form of an abstract power method, see \cite[Section 2.2]{KNY1987}.

These two estimates can be modified for the subspace iteration
implemented with implicit deflation and the exact shift-invert preconditioner
$K=(H-\sigma S)^{-1}$ for $\sigma<\la_i$.
Therein the current Ritz basis matrix $Z\in\CC^{n \times \widetilde{k}}$
is updated by $S$-orthonormal Ritz vectors associated with
the $i$th to the $(i{-}1{+}\widetilde{k})$th smallest Ritz values
in the subspace $\spanl\{U,(H-\sigma S)^{-1}SZ\}$. By using
a transformation with the representation $Z=VC$
and the $S$-orthogonal projector $Q=VV^*S$, we have
\[
Q(H-\sigma S)^{-1}SZ=VV^*S(H-\sigma S)^{-1}SVC
 \stackrel{\eqref{irf}}{=}VV^*S\big(V(\Lambda - \sigma \widetilde{I})^{-1}\big)C
 =V(\Lambda - \sigma \widetilde{I})^{-1}C.
\]
Consequently, we have
\[
\spanl\{U,(H-\sigma S)^{-1}SZ\}
 =\spanl\{U\}\oplus V \cdot \spanl\{(\Lambda - \sigma \widetilde{I})^{-1}C\}.
\] 
Thus the restricted iteration within $\spanl\{V\}$ reads
\begin{equation}\label{biid}
\spanl\{C'\}=\spanl\{(\Lambda - \sigma \widetilde{I})^{-1}C\},
\end{equation}
i.e., the block inverse iteration for the diagonal matrix
$(\Lambda - \sigma \widetilde{I})^{-1}$.
Then applying the analysis from \cite{PAR1980,KNY1987} gives the ratio
$(\la_{i-1+t}-\sigma)/(\la_{i+\widetilde{k}}-\sigma)$.  

Furthermore, the restricted form of the BPSD-id (Algorithm \ref{bpsdid})
has the trial subspace $\spanl\{C,\,(\Lambda - \sigma \widetilde{I})^{-1}C\}$
which is a superset of $\spanl\{C'\}$ from \eqref{biid}.
Thus Ritz value estimates for \eqref{biid} lead to indirect estimates for BPSD-id.

A more accurate estimate of this kind for BPSD-id can be shown
by using another auxiliary iteration whose restricted form reads
\begin{equation}\label{biida}
\spanl\{C'\}=\spanl\{(\Lambda - \sigma \widetilde{I})^{-1}C+\beta C\}.
\end{equation}
The accuracy can benefit from the choice of $\beta\in\RR$.

The following theorem predicts the cluster robustness of the BPSD-id. 

\begin{theorem}\label{bpsdidcr}
Consider the BPSD-id (Algorithm \ref{bpsdid})
with $K=(H-\sigma S)^{-1}$ and $\sigma<\la_i$. Denote by
$\theta_1^{(\ell)}\le\cdots\le\theta_{\widetilde{k}}^{(\ell)}$
the Ritz values in the $\ell$th subspace $\spanl\{Z^{(\ell)}\}$. Then
\begin{equation}\label{bpsdidcre}
 \frac{\theta_t^{(\ell)}-\la_{i-1+t}}{\la_n-\theta_t^{(\ell)}}
 \le\left(\frac{\kappa}{2-\kappa}\right)^{2\ell}
 \frac{\la_{i-1+t}-\sigma}{\la_n-\sigma}\,\tau
\end{equation}
with
\[\kappa=\left(\frac{\la_{i-1+t}-\sigma}{\la_{i+\widetilde{k}}-\sigma}\right)
 \left(\frac{\la_n-\la_{i+\widetilde{k}}}{\la_n-\la_{i-1+t}}\right)\]
and a constant $\tau>0$ depending on the initial subspace $\spanl\{Z^{(0)}\}$.
\end{theorem}
\begin{proof}
Following Lemma \ref{lmbpsdid1}, we observe the accompanying iteration
\begin{equation}\label{bpsdida}
\spanl\{C^{(\ell+1)}\}=\mathrm{RR}_{\Lambda_{\sigma},\widetilde{k},\uparrow}
\big(\spanl\{C^{(\ell)},\,\Lambda_{\sigma}^{-1}C^{(\ell)}\}\big)
\end{equation}
with the coefficient matrix $C^{(\ell)}=V^*SZ^{(\ell)}$
where the Rayleigh-Ritz procedure
$\mathrm{RR}_{\Lambda_{\sigma},\widetilde{k},\uparrow}(\cdot)$
extracts orthonormal Ritz vectors of
$\Lambda_{\sigma}=\mathrm{diag}(\lambda_i-\sigma,\ldots,\lambda_n-\sigma)$
associated with the $\widetilde{k}$ smallest Ritz values, i.e.,
$\theta_1^{(\ell+1)}{\,-\,}\sigma,\ldots,\theta_{\widetilde{k}}^{(\ell+1)}{\,-\,}\sigma$
are Ritz values of $\Lambda_{\sigma}$ in $\spanl\{C^{(\ell+1)}\}$, and
coincide with the $\widetilde{k}$ smallest Ritz values of $\Lambda_{\sigma}$
in $\spanl\{C^{(\ell)},\,\Lambda_{\sigma}^{-1}C^{(\ell)}\}$.

In addition, $\sigma<\la_i$ ensures that $\Lambda_{\sigma}$ is positive definite.
By using its square root matrix $\Lambda_{\sigma}^{1/2}$,
we define $E^{(\ell)}=\Lambda_{\sigma}^{1/2}C^{(\ell)}$
and $A=\Lambda_{\sigma}^{-1}$ so that
\[\spanl\{E^{(\ell)},\,AE^{(\ell)}\}
 =\spanl\{\Lambda_{\sigma}^{1/2}C^{(\ell)},\,\Lambda_{\sigma}^{-1/2}C^{(\ell)}\}
 =\Lambda_{\sigma}^{1/2}\spanl\{C^{(\ell)},\,\Lambda_{\sigma}^{-1}C^{(\ell)}\}.\]
Then the iteration \eqref{bpsdida} is equivalent to
\begin{equation}\label{bpsdidb}
\spanl\{E^{(\ell+1)}\}=\mathrm{RR}_{A,\widetilde{k},\downarrow}
\big(\spanl\{E^{(\ell)},\,AE^{(\ell)}\}\big)
\end{equation}
where the Rayleigh-Ritz procedure
$\mathrm{RR}_{A,\widetilde{k},\downarrow}(\cdot)$
extracts $A$-orthonormal Ritz vectors of $A$
associated with the $\widetilde{k}$ largest Ritz values.
For explaining this equivalence, we can denote by $C$ a basis matrix
of $\spanl\{C^{(\ell)},\,\Lambda_{\sigma}^{-1}C^{(\ell)}\}$
consisting of orthonormal Ritz vectors of $\Lambda_{\sigma}$.
Then $C^*C$ is an identity matrix,
and $C^*\Lambda_{\sigma}C$ is a diagonal matrix containing
Ritz values of $\Lambda_{\sigma}$. The corresponding
$E=\Lambda_{\sigma}^{1/2}C$ fulfills
\[E^*E=C^*\Lambda_{\sigma}C,\quad
 E^*AE=C^*C\]
so that its columns are $A$-orthonormal Ritz vectors of $A$
in $\spanl\{E\}=\Lambda_{\sigma}^{1/2}\spanl\{C\}$.
Moreover, the concerned Ritz values of $A$
are reciprocals of those of $\Lambda_{\sigma}$.

We further denote by $\alpha_1\ge\cdots\ge\alpha_{n-i+1}$
the eigenvalues of $A$, and by
$\psi_1^{(\ell)}\ge\cdots\ge\psi_{\widetilde{k}}^{(\ell)}$
the Ritz values in the $\ell$th subspace from \eqref{bpsdidb}. Then
\[\alpha_t=(\la_{i-1+t}-\sigma)^{-1},\quad
 \psi_t^{(\ell)}=\big(\theta_t^{(\ell)}-\sigma\big)^{-1}\]
so that the estimate \eqref{bpsdidcre} is equivalent to
\begin{equation}\label{bpsdidcre1}
 \frac{\alpha_t-\psi_t^{(\ell)}}{\psi_t^{(\ell)}-\alpha_{n-i+1}}
 \le\left(\frac{\kappa}{2-\kappa}\right)^{2\ell}\,\tau
\end{equation}
with
\[\kappa=\frac{\alpha_{\widetilde{k}+1}-\alpha_{n-i+1}}
 {\alpha_t-\alpha_{n-i+1}}.\]

Finally, we derive \eqref{bpsdidcre1} via the simplified version
\[
\spanl\{E^{(\ell+1)}\}=\spanl\{AE^{(\ell)}+\beta E^{(\ell)}\}
\]
of \eqref{bpsdidb}. Therein $\spanl\{E^{(\ell)}\}$ can be represented by
$f(A)\,\spanl\{E^{(0)}\}$ with $f(\eta)=\eta+\beta$.
This corresponds to the block form of an abstract power method
so that its convergence behavior
can be analyzed as in \cite[Section 2.2]{KNY1987}. By using
$\beta=-\tfrac12\big(\alpha_{t+1}+\alpha_{n-i+1}\big)$, we get
\[|f(\alpha_1)|\ge\cdots\ge|f(\alpha_t)|\ge|f(\alpha_j)|
 \quad\forall\ j\in\{t{\,+\,}1,\,\ldots,\,n{\,-\,}i{\,+\,}1\}.\]
Then the estimate \cite[(2.20)]{KNY1987} implies \eqref{bpsdidcre1}
where the constant $\tau$ is given by
the tangent square of the angle between
$\spanl\{E^{(0)}\}$ and the invariant subspace of $A$
associated with the $\widetilde{k}$ largest eigenvalues.
This angle depends on $\spanl\{Z^{(0)}\}$
since $E^{(0)}=\Lambda_{\sigma}^{1/2}V^*SZ^{(0)}$.
\end{proof}

Theorem \ref{bpsdidcr} 
predicts the cluster robustness of the BPSD-id 
since it indicates that the convergence rate
increases with the Ritz value index $t$ and 
the parameter $\kappa$ can be bounded away from $1$ for 
a sufficiently large block size $\widetilde{k}$. 

It is desirable to extend Theorem \ref{bpsdidcr}
to effectively positive definite preconditioners where
$\big(\kappa+\varepsilon(2-\kappa)\big)^{2\ell}
/\big((2-\kappa)+\varepsilon\kappa\big)^{2\ell}$
is a suitable convergence factor
with the quality parameter $\varepsilon$ defined in \eqref{pcde}.
Toward this aim, the analysis of the cluster robustness
of a preconditioned subspace iteration in \cite{NEZ2019}
can be adapted to the restricted formulation of the BPSD-id 
and yields an estimate of the form 
\begin{equation}\label{est2}
 \frac{\theta_t^{(\ell)}-\la_{i-1+t}}{\la_{i+\widetilde{k}}-\theta_t^{(\ell)}}
 \le\left(\varepsilon+(1-\varepsilon)\frac{\la_{i-1+t}-\nu}{\la_{i+\widetilde{k}}-\nu}\right)^{2\ell}
 \ \frac{\theta_{\widetilde{k}}^{(0)}-\la_{i-1+t}}
 {\la_{i+\widetilde{k}}-\theta_{\widetilde{k}}^{(0)}}
\end{equation}
for an arbitrary $\nu<\la_i$.
Therein the quality parameter $\varepsilon$ is related to $\nu$ and certain auxiliary vectors,
but close to that from an error propagation formulation such as \eqref{pcde}.

A recent approach in \cite{NEZ2022} for analyzing the cluster robustness
of the BPSD is also compatible with the restricted formulation of the BPSD-id.
The resulting estimate reads
\begin{equation}\label{est3}
 \frac{\theta_t^{(\ell)}-\la_{i-1+t}}{\la_{i+\widetilde{k}}-\theta_t^{(\ell)}}
 \le\left(\frac{\kappa+\varepsilon(2-\kappa)}{(2-\kappa)+\varepsilon\kappa}\right)^{2\ell}
 \ \frac{\theta_{\widetilde{k}}^{(0)}-\la_{i-1+t}}
 {\la_{i+\widetilde{k}}-\theta_{\widetilde{k}}^{(0)}}
\end{equation}
with
\[
\kappa=\left(\frac{\la_{i-1+t}-\nu}{\la_{i+\widetilde{k}}-\nu}\right)
 \left(\frac{\la_n-\la_{i+\widetilde{k}}}{\la_n-\la_{i-1+t}}\right).
\]
In the case $\varepsilon=0$ and $\nu=\sigma$,
the convergence factor in \eqref{est3} coincides with that in \eqref{bpsdidcre}.

An alternative of the multi-step estimates \eqref{est2} and \eqref{est3}
follows from the single-step estimate \eqref{eq:bpsdidbd}
in Theorem~\ref{bpsdidth} for $j=i{\,-\,}1{\,+\,}t$, namely,
\begin{equation}\label{est1}
 \frac{\theta_t^{(\ell)}-\la_{i-1+t}}{\la_{i+t}-\theta_t^{(\ell)}}
 \le\left(\frac{\kappa+\varepsilon(2-\kappa)}{(2-\kappa)+\varepsilon\kappa}\right)^{2\ell}
 \frac{\theta_t^{(0)}-\la_{i-1+t}}{\la_{i+t}-\theta_t^{(0)}}
\end{equation}
with
\[\kappa=\left(\frac{\la_{i-1+t}-\nu}{\la_{i+t}-\nu}\right)
 \left(\frac{\la_n-\la_{i+t}}{\la_n-\la_{i-1+t}}\right).\]
The estimate \eqref{est1} cannot predict the cluster robustness,
but can provide better bounds for the first steps
in comparison to \eqref{est2} and \eqref{est3}.

\section{Numerical experiments} \label{sec:numerics} 

In this section, we illustrate main convergence estimates
of the PSD-id and the BPSD-id presented in the previous sections by
several numerical examples. Example~\ref{exe1} with non-clustered
targets eigenvalues demonstrates single-step estimates
for the BPSD-id in Section~\ref{sec:bpsdid}.
In Example~\ref{exe4}, we implement the BPSD-id
for a large-scale eigenvalue problem and compare
three effectively positive definite preconditioners.
In Example~\ref{exe3}, we revisit
\cite[Example 5.2]{CBPS2018} where the eigenvalue problem
is derived by the partition-of-unity finite element
for a self-consistent pseu-dopotential density functional calculation.
We refine the implementation of the PSD-id therein (by accelerating inner steps)
and interpret a cubic convergence for dynamic larger shifts
by estimates in Section~\ref{sec:sigma}.
Example~\ref{exe2} with clustered targets eigenvalues
demonstrates multi-step estimates for the BPSD-id
in Section~\ref{sec:multi}.

\begin{example}\label{exe1}
We consider the Laplacian eigenvalue problem 
$-\Delta u = \la u$ on the rectangle
$[0,\,1.5]\times[0,\,1]$ with two slits $\{0.5\}\times[0.45,\,0.55]$ 
and $\{1\}\times[0.45,\,0.55]$; see Figure~\ref{fig1}.
The boundary condition
on the rectangle boundary and the slits is simply $u=0$. 
By using the five-point star discretization with the mesh size $h=45/3600$,
we get the eigenvalue problem \eqref{evp} of order $n=9383$ and $S=I$.
The seven smallest eigenvalues of $(H,S)$ are well separated:
\[\la_1\approx27.07834,\quad\la_2\approx38.24327,\quad\la_3\approx45.24858,\]
\[\la_4\approx49.32646,\quad\la_5\approx58.36810,\quad\la_6\approx78.91626,
 \quad\la_7\approx89.70648.\]

\begin{figure}[htbp]
\quad\vspace{1em}
\begin{center}
\quad\includegraphics[width=\textwidth]{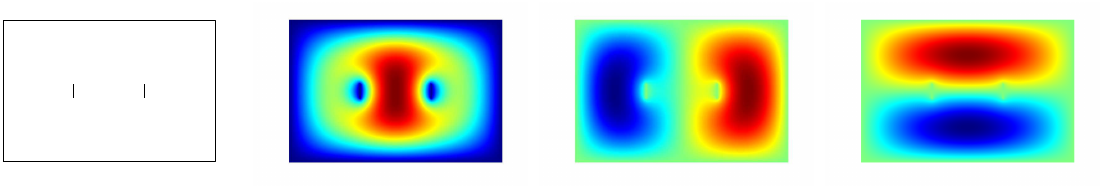}
\end{center}
\vspace{1em}
\caption{\small
Laplacian eigenvalue problem in Example~\ref{exe1}.
The two slits on the domain clearly influence the shapes of eigenfunctions
associated with the three smallest operator eigenvalues.}
\label{fig1}
\end{figure}

We compute the six smallest eigenvalues 
by the BPSD-id (Algorithm \ref{bpsdid}) with 
$\{k,\widetilde{k}\} = \{1,2\}, \{2,3\}, \{3,4\}$.
Recall that $k$ is the number of wanted eigenvalues in each run
and $\widetilde{k}$ is the block size. 
For instance, if $\{k,\widetilde{k}\} = \{2,3\}$, then 
the BPSD-id with block size 3 computes 2 eigenvalues
in each of three successive runs.

We first consider fixed shifts
and use incomplete matrix factorizations
for generating the preconditioner $K$, namely, 
\begin{center}
\verb|ichol(H-sigma*S,struct('type','ict','droptol',3e-5))|
\end{center}
with the shift $\sigma=20$ for the first run (i.e. for $i=1$) and
\begin{center}
\verb|ilu(H-sigma*S,struct('type','crout','milu','row','droptol',3e-5))|
\end{center}
with the shift $\sigma=\la_{i-1}$ for further runs.

Figure~\ref{fig1a} depicts the convergence
of the residual norm $\|R(:,\,1\,{:}\,k)\|_{S^{-1}}$.
We run $1000$ random initial subspaces and show the slowest convergence.
We observe that the residual norm decreases monotonically
for different choices of $\{k,\widetilde{k}\}$.

\begin{figure}[htbp]
\begin{center}
\includegraphics[width=\textwidth]{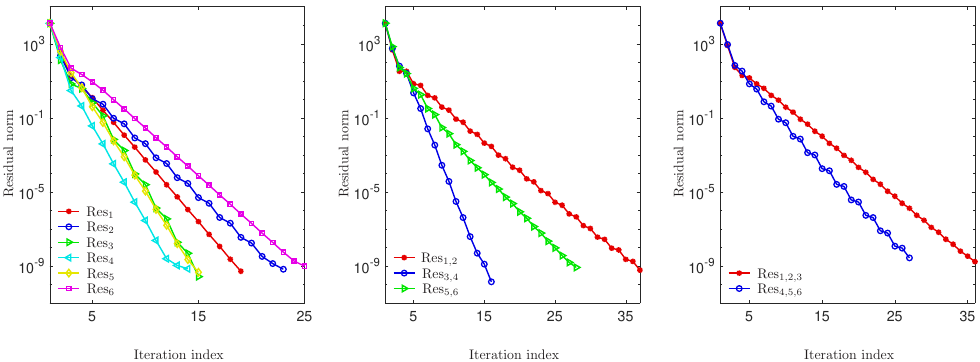}\quad
\end{center}
\par\vskip -2ex
\caption{\small Convergence of the BPSD-id (Algorithm \ref{bpsdid})
with respect to the residual norm $\|R(:,\,1\,{:}\,k)\|_{S^{-1}}$
for computing the six smallest eigenvalues in Example~\ref{exe1}.
The pair $\{k,\widetilde{k}\}$ is $\{1,2\}$ (left),
$\{2,3\}$ (center) and $\{3,4\}$ (right).}
\label{fig1a}
\end{figure}

The convergence of the BPSD-id with respect to
the Ritz value error $\theta_t-\la_t$ for $t\in\{1,\ldots,6\}$ 
is depicted in Figure~\ref{fig1b},
where the indices of Ritz values are permutated
to match the target eigenvalues. 
We observe that Theorem \ref{bpsdidth} provides 
sharp error bounds in dotted curves.  
For evaluating these bounds, we determine
the quality parameter $\varepsilon$
by using Lemma \ref{pcd} and the fact that
the nonzero eigenvalues of 
$\widetilde{K}\Lambda_{\sigma}=(V^*SKSV)V^*(H-\sigma S)V$
coincide with those of $K(SVV^*)(H-\sigma S)VV^*S$.
Then the extremal eigenvalues of $\widetilde{K}\Lambda_{\sigma}$
can be obtained by {\tt eigs} applied to a subroutine 
for matrix-vector multiplications
where $(VV^*S)x$ is computed by \,$x-U\big(U^*(Sx)\big)$,
and $(SVV^*)y$ by \,$y-S\big(U(U^*y)\big)$.

We note that the sharpness statement for
the estimate \eqref{eq:bpsdidbd} cannot easily be observed
for random iterates as shown in Figure~\ref{fig1b}.
The attainability of the associated equality
within certain low-dimensional invariant subspaces
can however be illustrated similarly to \cite[Figure 4.5]{NEZ2014}.

\begin{figure}[htbp]
\begin{center}
\includegraphics[width=\textwidth]{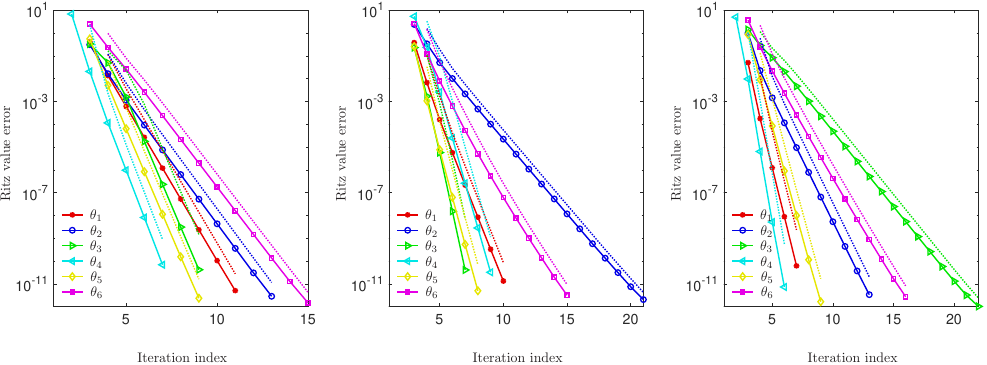}\quad
\end{center}
\par\vskip -2ex
\caption{\small Convergence of the BPSD-id (Algorithm \ref{bpsdid})
with respect to the Ritz value error $\theta_t-\la_t$
for computing the six smallest eigenvalues in Example~\ref{exe1}
with fixed shifts.
The error bounds in dotted curves are determined
by Theorem \ref{bpsdidth}.}
\label{fig1b}
\end{figure}

Let us consider refining the preconditioner $K$ by a dynamic shift $\sigma$.
With the index $i$ of the smallest target eigenvalue
in the current run, we estimate the ratio 
$\eta=(\theta_i-\la_i)/(\la_{i+1}-\theta_i)$
roughly by $(\theta_{i,\mathrm{old}}-\theta_i)/(\theta_{i+1}-\theta_i)$
as suggested for the PSD-id in \cite[Section 5]{CBPS2018}.
If $\eta$ and the residual norm $\|R(:,\,1\,{:}\,k)\|_{S^{-1}}$
are both smaller than the threshold $0.1$, we update $\sigma$ by
$\sigma \gets (\sigma+\theta_i)/2$ and refine the preconditioner $K$ by 
\begin{center}
\verb|ilu(H-sigma*S,struct('type','crout','milu','row','droptol',max(eta,1e-12)))|\ .
\end{center}
We mark the first refinement by ``+''. As shown in Figure \ref{fig1c},
this modification leads to an acceleration.
The improvement is evident for the $\{k,\widetilde{k}\}=\{1,2\}$
(left subfigure). Furthermore, Theorem~\ref{bpsdidth}
produces proper bounds in dotted curves.

\begin{figure}[htbp]
\begin{center}
\includegraphics[width=\textwidth]{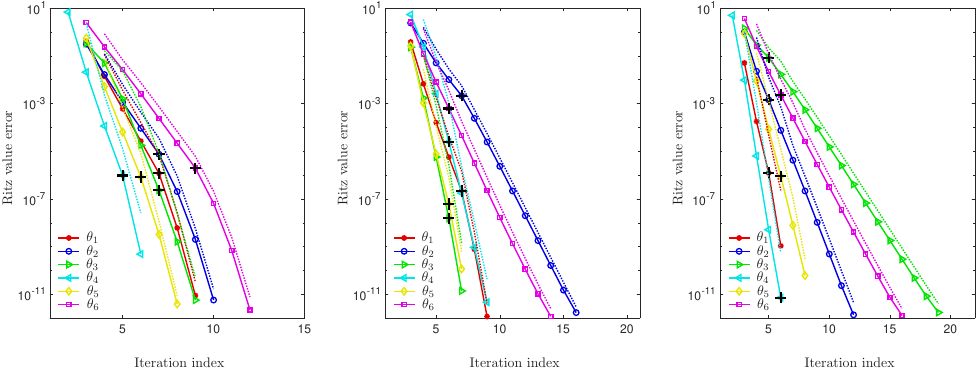}\quad
\end{center}
\par\vskip -2ex
\caption{\small 
Convergence of the BPSD-id (Algorithm \ref{bpsdid})
with respect to the Ritz value error $\theta_t-\la_t$
for computing the six smallest eigenvalues in Example~\ref{exe1}
with dynamic shifts.}
\label{fig1c}
\end{figure}

\end{example}

\begin{example}\label{exe4}

For discussing the performance of the BPSD-id in large-scale problems,
let us consider a matrix pair arising from an adaptive finite element discretization
of the Laplacian eigenvalue problem on a wrench-shaped domain
with homogeneous Dirichlet boundary conditions; see Figure~\ref{fig4a}.
The boundary is defined by
\[
 \big\{20\cos(t)+10\cos(2t),\ 5\sin(t)+\sin(5t);\ t\in[0,\,2\pi)\big\}.
\]
Similarly to \cite[Appendix]{NEZ2019}, matrix eigenvalue problems
are generated successively by an adaptive finite element discretization.
The refinement is controlled by residuals of approximate eigenfunctions
associated with the three smallest operator eigenvalues.

\begin{figure}[htbp]
\quad\vspace{1em}
\begin{center}
\includegraphics[width=\textwidth]{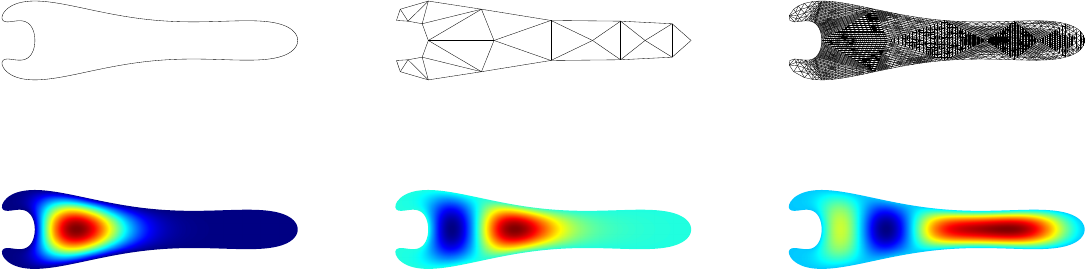}\quad
\end{center}
\vspace{1em}
\caption{\small
Laplacian eigenvalue problem in Example~\ref{exe4}.
The first row displays the domain, the initial grid and 
an adaptively refined grid.
The second row illustrates approximate eigenfunctions
associated with the three smallest operator eigenvalues.
The corresponding residuals are used for the grid refinement.}
\label{fig4a}
\end{figure}

We consider the matrix pair $(H,S)$ from the $29$th grid of the discretization
with $n=1{,}618{,}797$ degrees of freedom.
The seven smallest eigenvalues of $(H,S)$ are well separated
and located in the interval $(0.1418824,\,0.4562653)$.
Similarly to Example~\ref{exe1}, we compute the six smallest eigenvalues
$\la_1,\ldots,\la_6$ by $3$ successive runs of the 
BPSD-id (Algorithm \ref{bpsdid}) with $\{k,\widetilde{k}\}=\{2,4\}$.

Figure~\ref{fig4b} illustrates the convergence of the BPSD-id with respect to
the Ritz value error $\theta_t-\la_t$ for $t\in\{1,\ldots,6\}$.
Therein we compare three effectively positive definite preconditioners.

In the left subfigure, the preconditioner $K$ is generated by
the incomplete matrix factorization
\begin{center}
\verb|ichol(H-sigma*S,struct('type','ict','droptol',1e-6))|
\end{center}
with the shift $\sigma=0.1<\la_1$. We use this $K$ for each run,
since generating $K$ by \verb|ilu| with $\sigma=\la_{i-1}$ is too costly.
Consequently, more outer steps are required in the second and third runs
(curves for $\theta_3,\theta_4$ and $\theta_5,\theta_6$)
than in the first run (curves for $\theta_1,\theta_2$).

In the central subfigure, we modify $K$ for the second and third runs
where approximations of $(H-\sigma S)^{-1}$ for $\sigma=\la_{i-1}$
are constructed by MINRES with the above \verb|ichol| factorization
and the tolerance $0.1$.  This clearly reduces the number of required 
outer steps.  However, the computational time for these two runs
reads 176 seconds, longer than 105 seconds measured in the left subfigure.

The same modification with the tolerance $0.01$, as illustrated in the right subfigure,
does not lead to a further significant reduction of step numbers,
but increases the computational time to 208 seconds.
In addition, modifications with dynamic shifts
are also less efficient with respect to the total computational time.
Thus a simply applicable preconditioner with fixed shifts is occasionally
more appropriate.

\begin{figure}[htbp]
\begin{center}
\includegraphics[width=\textwidth]{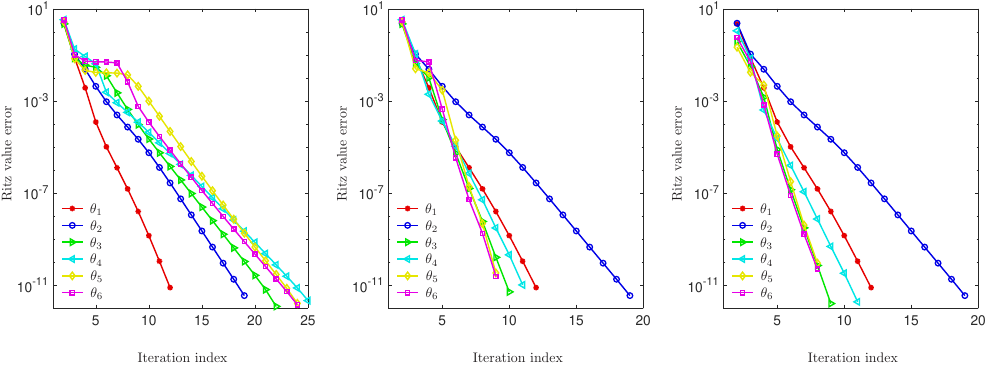}\quad
\end{center}
\par\vskip -2ex
\caption{\small Convergence of the BPSD-id (Algorithm \ref{bpsdid})
with respect to the Ritz value error $\theta_t-\la_t$
for computing the six smallest eigenvalues in Example \ref{exe4}.
Left: Using \texttt{ichol} for generating the preconditioner with a constant shift.
Center: Modifying the preconditioner by MINRES with the tolerance $0.1$.
Right: Modifying the preconditioner by MINRES with the tolerance $0.01$.}
\label{fig4b}
\end{figure}
\end{example}

\begin{example}\label{exe3}

To verify the sharpness of the estimates with larger shifts 
discussed in Section~\ref{sec:sigma}, we use the matrix pencil $(H,S)$ 
of order $n=5336$ derived from partition-of-unity finite element method
for quantum-mechanical materials calculation in \cite[Example 5.2]{CBPS2018}.
The matrix $H$ is given by a rank-$l$ modification of an $n \times n$
sparse matrix $\widehat{H}$ with $l=26$.
Both $H$ and $S$ are ill-conditioned and their condition numbers
are $\mathcal{O}(10^{10})$. Furthermore, $H$ and $S$ share a common 
near-nullspace $\mbox{span}\{V\}$ of dimension 1000
such that $\|HV\| = \|SV\| = \mathcal{O}(10^{-4})$. This is considered 
as an extremely ill-conditioned eigenvalue problem.  

For the BPSD-id, the explicit form of $H$
is dense so that constructing $K$ by incomplete matrix factorizations
is not efficient. 
Thus MINRES with preconditioner \verb|chol(S)| is used for computing
the preconditioned residual $Kr$ with $K\approx(H-\sigma S)^{-1}$.
The stopping criterion of MINRES uses the residual norm
$\psi=\|r\|_2/(\|Hz\|_2+\|\rho(z)Sz\|_2)$ instead of $\|r\|_{S^{-1}}$.
The initial $\sigma$ is $-1$, which   
is smaller than the smallest eigenvalue $\la_1\approx-0.8888$.
The tolerance of MINRES is $0.1$ for $\sigma = -1$.
If $\psi$ and an estimated value of 
$\big(\rho(z)-\la_i\big)/\big(\la_{i+1}-\rho(z)\big)$ are sufficiently small, 
then $\sigma$ is set equal to $\rho(z)$.

We begin with the PSD-id (Algorithm \ref{psdid}) implemented as a 
weakened form of the BPSD-id where $r=R(:,1)$ is used for
constructing the trial subspace.  This corresponds to an acceleration 
of the PSD-id so that Theorem~\ref{psdidrq} is still applicable
and the estimate \eqref{psdidrqe} indicates a cubic convergence
for larger shifts from the interval $\big(\la_i,\,(\la_i+\la_{i+1})/2\big)$.

We note that the inner steps (MINRES) in the above implementation
can significantly be accelerated by using $\widehat{H}+3S$ instead of
$S$ in \verb|chol(S)|. In addition, a substantial acceleration
of the outer steps in the third and fourth runs is enabled
by using $\la_{i-1}$ as the initial $\sigma$.  This improvement can 
be observed by comparing the reduction of the residual norm $\psi$
in Figure~\ref{fig2a} (left) with that in \cite[Figure 5.2]{CBPS2018}.

\begin{figure}[htbp]
\begin{center}
\includegraphics[width=\textwidth]{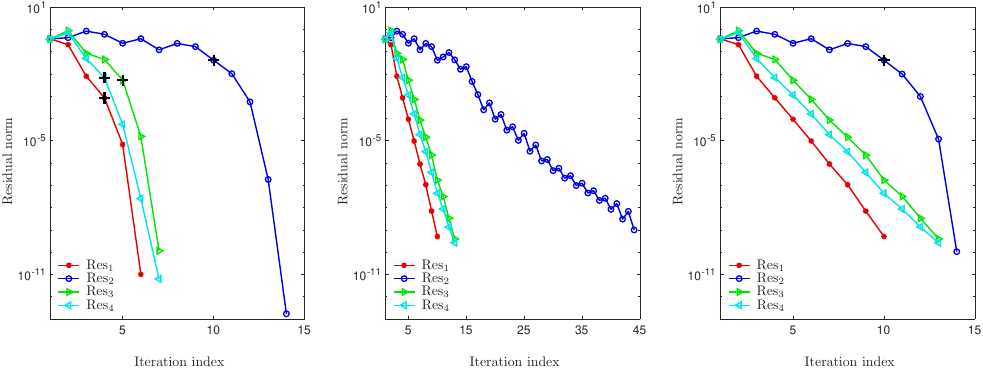}\quad
\end{center}
\par\vskip -2ex
\caption{\small Convergence of a weakened form of BPSD-id
using only one residual vector and an relative residual norm
for computing the four smallest eigenvalues in Example \ref{exe3}.
The preconditioner is constructed by MINRES with certain shifts.
Left: Enabling dynamic shifts. The first dynamic step is marked by ``+''.
Center: Using fixed shifts only. Right: Enabling dynamic shifts
when the convergence in the first steps is slow.}
\label{fig2a}
\end{figure}

Furthermore, Figure~\ref{fig2a} (center) depicts the reduction
of $\psi$ for fixed shifts, i.e., $\sigma=-1$ for the first run
and $\sigma=\la_{i-1}$ for further runs. The tolerance of MINRES
is constantly $0.1$. In comparison to the version with dynamic shifts,
only the second run is considerably slowed down with respect to
the outer steps. The total computational time is however only slightly increased.
In Figure \ref{fig2a} (right), a hybrid version is implemented where
dynamic shifts are used if the iteration index for the first possible switch
is larger than $6$. Then dynamic shifts are enabled
in the second run, and the total computational time
is reduced in comparison to the previous two versions.

Since our estimates are formulated for Ritz value errors,
we observe in Figure~\ref{fig2b} the error $\theta_t-\la_t$
for the three versions mentioned above. Therein $\theta_t$
denotes the approximate eigenvalue in the $t$\,th run.
The convergence is obviously (piecewise) linear for fixed shifts
and can become cubic by using dynamic shifts, as predicted
by Theorem \ref{psdidth} and Theorem \ref{psdidrq}, respectively.

\begin{figure}[htbp]
\begin{center}
\includegraphics[width=\textwidth]{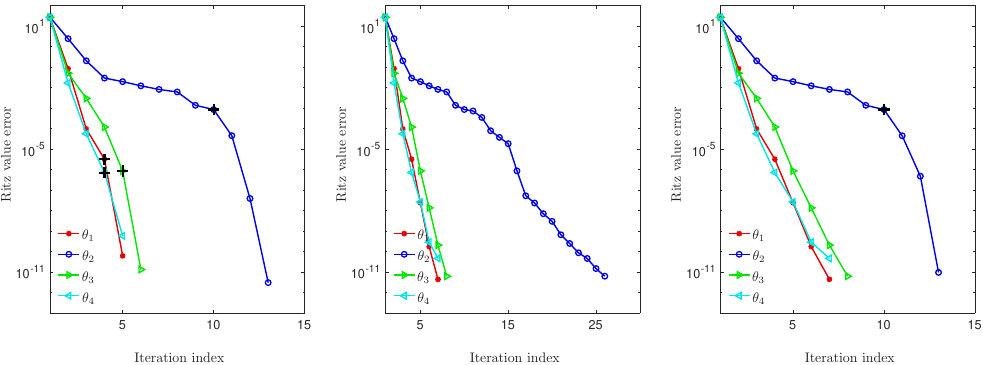}\quad
\end{center}
\par\vskip -2ex
\caption{\small Convergence of a weakened form of BPSD-id
with respect to the Ritz value error $\theta_t-\la_t$
in addition to Figure \ref{fig2a}.}
\label{fig2b}
\end{figure}

When implementing the standard BPSD-id with shifts
from $\big(\la_i,\,(\la_i+\la_{i+1})/2\big)$ in Example~\ref{exe3},
a cubic convergence can occur for the smallest Ritz value in each run.
Further Ritz values only converge linearly. This difference matches
the two estimates from Theorem~\ref{bpsdidrq}.
\end{example}

\begin{example}\label{exe2}

By using a block size which is larger than the cluster size,
the BPSD-id can efficiently compute clustered eigenvalues.
This fact has been analyzed in Theorem~\ref{bpsdidcr}
for exact shift-invert preconditioning. The estimate \eqref{bpsdidcre}
corresponds to a special form of \eqref{est3} concerning
effectively positive definite preconditioners. The estimate \eqref{est3}
can be derived by adapting the analysis from \cite{NEZ2022}
to the restricted formulation of the BPSD-id analogously to
Section \ref{sec:bpsdid}. In this example, we compare \eqref{est3}
with its counterparts \eqref{est2} and \eqref{est1}
which are based on \cite{NEZ2019}
and the single-step estimate \eqref{eq:bpsdidbd}, respectively.

We modify the eigenvalue problem in Example \ref{exe1} by setting larger slits
$\{0.5\}\times[0.1,\,0.9]$ and $\{1\}\times[0.1,\,0.9]$
on the rectangle $[0,\,1.5]\times[0,\,1]$; see Figure~\ref{fig2}.
The mesh size $h=1/80$ leads to $n=9271$.
The modified domain can be regarded
as three small rectangles connected by narrow gates. For a sufficiently small
gate width, the eigenvalue problem is almost split into three partial problems
of the same size. Thus the eigenvalues
are roughly copies of those from partial problems.
As a result, there are two tight clusters among the seven smallest eigenvalues:
\[
 \la_1,\la_2,\la_3\in(49.24886,\,49.32647),\quad
 \la_4,\la_5,\la_6\in(78.61283,\,78.91626),\quad
 \la_7\approx127.5209.
\]

\begin{figure}[htbp]
\quad\vspace{1em}
\begin{center}
\includegraphics[width=\textwidth]{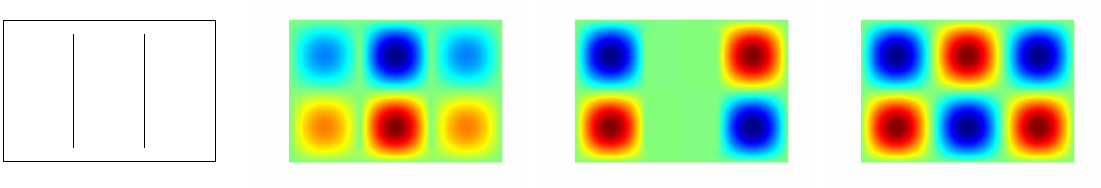}
\end{center}
\vspace{1em}
\caption{\small
Laplacian eigenvalue problem in Example~\ref{exe2}.
The eigenfunctions associated with the fourth to sixth
smallest operator eigenvalues are displayed.}
\label{fig2}
\end{figure}

Let us examine a run of the BPSD-id with $i=4$ and $\widetilde{k}=3$, 
which means that the eigenvalues
$\la_4,\la_5,\la_6$ are to be computed. The preconditioner $K$ is generated by
\begin{center}
\verb|ilu(H-sigma*S,struct('type','crout','milu','row','droptol',3e-5))|
\end{center}
for $\sigma=\la_3$. We denote by $\mathrm{Bound}_1$,
$\mathrm{Bound}_2$ and $\mathrm{Bound}_3$ the bounds
for the Ritz value error $\theta_t-\la_t$ (with simplified indices)
which are determined on the basis of
\eqref{est1}, \eqref{est2} and \eqref{est3}, respectively.
We compare these bounds with the error $\theta_t-\la_t$ in Figure \ref{fig3}
with subfigures for $t\in\{4,\,5,\,6\}$ by evaluating their numerical maxima
concerning $1000$ random initial subspaces.

\begin{figure}[htbp]
\begin{center}
\includegraphics[width=\textwidth]{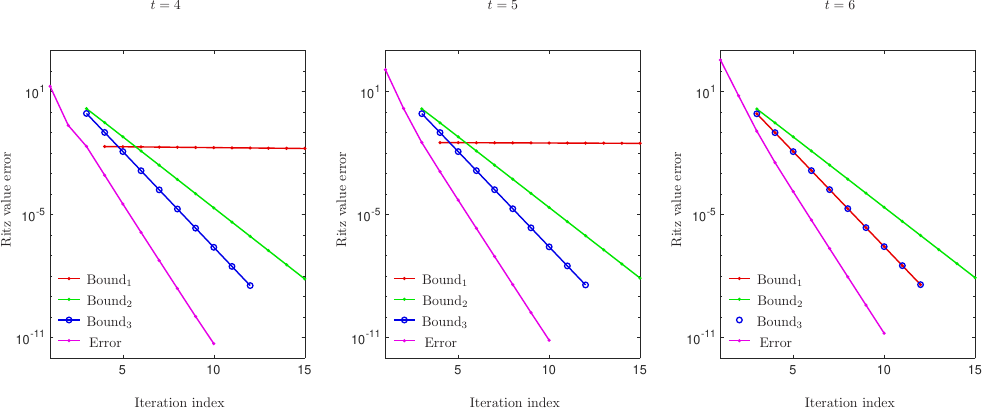}\quad
\end{center}
\par\vskip -2ex
\caption{\small Convergence of the BPSD-id (Algorithm \ref{bpsdid})
with respect to the Ritz value error $\theta_t-\la_t$ for $t\in\{4,\,5,\,6\}$
in Example \ref{exe2} with clustered eigenvalues.
The three bounds are based on the multi-step estimates
\eqref{est1}, \eqref{est2} and \eqref{est3}.}
\label{fig3}
\end{figure}

We note that  $\mathrm{Bound}_2$ and $\mathrm{Bound}_3$ are nearly invariant
for $t$ due to the eigenvalue cluster $\{\la_4,\la_5,\la_6\}$.
Moreover, $\mathrm{Bound}_3$ is evidently more accurate than $\mathrm{Bound}_2$.
For $t\in\{4,\,5\}$, the clustered eigenvalues make $\mathrm{Bound}_1$ nearly constant.
For $t=6$, the moderate gap between $\la_6$ and $\la_7$ leads to a meaningful
$\mathrm{Bound}_1$ which actually coincides with $\mathrm{Bound}_3$.

\end{example}

\section{Conclusion}

The limitation of the approach presented in \cite{CBPS2018}
for analyzing the convergence behavior of the PSD-id method is overcome by
embedding concise bounds from \cite{NEY2012,NEZ2014} concerning
the PSD and BPSD methods. The new estimates are more flexible
with weaker assumptions, natural description of preconditioning
and extension to block iterations. Therein the preconditioners
are assumed to be effectively positive definite and particularly include
approximative shift-invert preconditioners where the shift is smaller than
the target eigenvalue. In addition, the case of larger shifts is discussed
based on the analysis of an abstract power method from \cite{KNY1987}
and the analysis of an inexact Rayleigh quotient iteration from \cite{NTY2003}. 
Furthermore, the cluster robustness of the BPSD-id
is analyzed for exact shift-invert preconditioning
analogously to an abstract block iteration from \cite{KNY1987}.
A more general analysis for effectively positive definite preconditioners
is enabled by recent progress for the BPSD from \cite{NEZ2022}.
Topics for further study include ``implicit deflation'' versions
of the LOBPCG and various Davidson methods
as well as practical settings of shifts and block sizes.

\bibliographystyle{plain}

\begin{thebibliography}{10}

\itemsep1ex

\bibitem{BPK1996}
J.H.~Bramble, J.E.~Pasciak, and A.V.~Knyazev,
\textit{A subspace preconditioning algorithm for 
eigenvector/eigenvalue computation},
Adv. Comput. Math.
{6} (1996), 159--189.

\bibitem{CBPS2018}
Y.~Cai, Z.~Bai, J.E.~Pask, and N.~Sukumar,
\textit{Convergence analysis of a locally accelerated preconditioned steepest
        descent method for Hermitian-definite generalized eigenvalue problems},
J. Comp. Math.
{36} (2018), 739--760.

\bibitem{KNY1987}
A.V.~Knyazev,
\textit{Convergence rate estimates for iterative methods for a mesh symmetric
        eigenvalue problem},
Russian J. Numer. Anal. Math. Modelling
{2} (1987), 371--396.

\bibitem{KNN2003e}
A.V.~Knyazev and K.~Neymeyr,
\textit{Efficient solution of symmetric eigenvalue problems using multigrid
preconditioners in the locally optimal block conjugate gradient method},
Electron. Trans. Numer. Anal. 
{15} (2003), 38--55.

\bibitem{KNS1991}
A.V.~Knyazev and A.L.~Skorokhodov,
\textit{On exact estimates of the convergence rate of the steepest ascent
  method in the symmetric eigenvalue problem},
Linear Algebra Appl.
{154--156} (1991), 245--257.

\bibitem{NEY2012}
K.~Neymeyr,
\textit{A geometric convergence theory for the preconditioned steepest
        descent iteration},
SIAM J. Numer. Anal.
{50} (2012), 3188--3207.

\bibitem{NOZ2011}
K.~Neymeyr, E.E.~Ovtchinnikov, and M.~Zhou,
\textit{Convergence analysis of gradient iterations for the symmetric
        eigenvalue problem},
SIAM J. Matrix Anal. Appl.
{32} (2011), 443--456.

\bibitem{NEZ2014}
K.~Neymeyr and M.~Zhou,
\textit{The block preconditioned steepest descent iteration for
        elliptic operator eigenvalue problems},
Electron. Trans. Numer. Anal.
{41} (2014), 93--108.

\bibitem{NTY2002}
Y.~Notay,
\textit{Combination of Jacobi-Davidson and conjugate gradients
        for the partial symmetric eigenproblem},
Numer. Linear Algebra Appl.
{9} (2002), 21--44.

\bibitem{NTY2003}
Y.~Notay,
\textit{Convergence analysis of inexact Rayleigh quotient iteration},
SIAM J. Matrix Anal. Appl.
{24} (2003), 627--644.

\bibitem{OVT2006}
E.E.~Ovtchinnikov,
\textit{Cluster robustness of preconditioned gradient
        subspace iteration eigensolvers},
Linear Algebra Appl.
{415} (2006), 140--166.

\bibitem{OVT2006s}
E.E.~Ovtchinnikov,
\textit{Sharp convergence estimates for the preconditioned steepest descent
  method for Hermitian eigenvalue problems},
SIAM J. Numer. Anal.
{43} (2006), 2668--2689.

\bibitem{PAR1980}
B.N.~Parlett,
\textit{The Symmetric Eigenvalue Problem},
Prentice-Hall, Englewood Cliffs, NJ, 1980.
Reprinted as Classics in Applied Mathematics 20, SIAM,
Philadelphia, 1997.

\bibitem{SAA1992}
Y.~Saad,
\textit{Numerical Methods for Large Eigenvalue Problems},
Manchester University Press, 1992.

\bibitem{STS1990}
G.~W.~Stewart and J.~Sun, 
\textit{Matrix Perturbation Theory},
Academic Press, 1990.

\bibitem{SAM1958}
B.A.~Samokish,
\textit{The steepest descent method for an eigenvalue problem with
  semi-bounded operators},
Izv. Vyssh. Uchebn. Zaved. Mat.
{5} (1958), 105--114 (in Russian).

\bibitem{NEZ2019}
M.~Zhou and K.~Neymeyr,
\textit{Cluster robust estimates for block gradient-type eigensolvers},
Math. Comp. 
{88} (2019), 2737--2765.

\bibitem{NEZ2022}
M.~Zhou and K.~Neymeyr,
\textit{Convergence rates of individual Ritz values
in block preconditioned gradient-type eigensolvers},
Technical Report \texttt{https://arxiv.org} \texttt{/abs/2206.00585}, 2022.

\end{thebibliography}

\end{document}